\documentclass[final]{amsart}
\usepackage{pdfsync}
\usepackage{amsmath,amssymb}
\usepackage{enumerate}
\usepackage{xspace}
\usepackage{boxedminipage}
\usepackage{algorithmicx}
\usepackage{algpseudocode}
\usepackage{listings}
\usepackage{tabularx}
\usepackage{xcolor}
\usepackage[figurename=Fig.,labelfont=bf,labelsep=period]{caption}
\usepackage{subcaption}
\usepackage{fullpage}
\usepackage[dvipsnames]{xcolor}
\usepackage[
  colorlinks=true,
  linkcolor=Purple,
  citecolor=Orange,
  urlcolor=Purple
]{hyperref}

\usepackage{tcolorbox}
\usepackage{smartdiagram}
\usepackage{colortbl}
\tcbuselibrary{skins}

\usepackage[backend=biber,style=numeric,giveninits, maxbibnames=9]{biblatex}
\tcbset{tab2/.style={enhanced,fonttitle=\bfseries,fontupper=\normalsize\sffamily,
colback=yellow!10!white,colframe=red!50!black,colbacktitle=red!40!white,
coltitle=black,center title}}

\usepackage{tikz,pgfplots} 
\usetikzlibrary{patterns}
\usetikzlibrary{spy}
\usetikzlibrary{cd}
\usepackage{pgfplots}
\pgfplotsset{width=5cm,compat=1.14}
\usepgfplotslibrary{colorbrewer}
\usetikzlibrary{intersections}
\usetikzlibrary{arrows,snakes,shapes}
\usetikzlibrary{decorations.markings}
\tikzset{
    set arrow inside/.code={\pgfqkeys{/tikz/arrow inside}{#1}},
    set arrow inside={end/.initial=>, opt/.initial=},
    /pgf/decoration/Mark/.style={
        mark/.expanded=at position #1 with
        {
            \noexpand\arrow[\pgfkeysvalueof{/tikz/arrow inside/opt}]{\pgfkeysvalueof{/tikz/arrow inside/end}}
        }
    },
    arrow inside/.style 2 args={
        set arrow inside={#1},
        postaction={
            decorate,decoration={
                markings,Mark/.list={#2}
            }
        }
    },
  }
  
\newcommand{\imagewithcentreline}[2][]{%
  \begin{tikzpicture}
    \node[anchor=south west, inner sep=0] (img) {\includegraphics[#1]{#2}};
    \begin{scope}[x={(img.south east)}, y={(img.north west)}]
      \draw[gray] (0.5, 0) -- (0.5, 1);
      \draw[gray] (0, 0.67) -- (1, 0.67);
    \end{scope}
  \end{tikzpicture}%
}
\definecolor{pone}{HTML}{2bc2c2}
\definecolor{ptwo}{HTML}{f4b811}
\definecolor{pthree}{HTML}{de663e}
\definecolor{pfour}{HTML}{4c90ba}

\definecolor{small-Wi}{HTML}{febd59}
\definecolor{large-Wi}{HTML}{00008b}

\newcommand{\ConvergenceTriangle}[4]{%
    \coordinate (A) at (axis cs:#1, #2); 
    \coordinate (B) at (axis cs:{#1*#3}, #2); 
    \coordinate (C) at (axis cs:{#1*#3}, {#2/#3^#4)}); 

    \draw[thick] (A) -- (B);
    \draw[thick] (A) -- (C);
    \draw[thick] (B) -- (C) node[midway, right] {#4};  
    
}

\usepackage{tristan}

\newcommand{\inflow}{\Gamma_-}
\newcommand{\outflow}{\Gamma_+}
\renewcommand{\Re}{\text{Re\xspace}}
\renewcommand{\De}{\text{De\xspace}}

\newcommand{\Wi}{\text{Wi\xspace}}

\newcommand{\sobhcurl}[1]{\ensuremath{\HH\qp{\curl;#1}}}
\newcommand{\sobhdiv}[1]{\ensuremath{\HH\qp{\div;#1}}}

\newcommand{\normal}{\vec n}
\renewcommand{\sig}{\vec \sigma}
\newcommand{\symgrad}[1]{\ensuremath{\vec \varepsilon(#1)}}
\renewcommand{\ddt}{\dfrac{\d}{\d t}}
\renewcommand{\avg}[1]{\mathinner{\doublelbrace\!\!#1\!\doublerbrace}}
\newcommand{\half}{\frac{1}{2}}
\newcommand{\ujump}[1]{\ensuremath{\left\lfloor\!\!\left\lfloor #1 \right\rfloor\!\!\right\rfloor}}

\numberwithin{equation}{section}
\setboolean{showtodo}{true}
\setboolean{shownotes}{true}
\setboolean{showchanges}{false}
\setboolean{usemathrsfs}{false}
\author{
  Ben S. Ashby
}
\address{
  Ben S. Ashby
  \thanks{
    Institute for Mathematical Innovation,
    University of Bath, Bath BA2 7AY, UK
    {\tt{bsa34@bath.ac.uk}}.
}}

\author{
  Gabriel R.  Barrenechea
  }
  \address{
  Gabriel R.  Barrenechea
  \thanks{
    Department of Mathematics and Statistics,
    University of Strathclyde, 26 Richmond Street, Glasgow G1 1XK, UK,
    {\tt{gabriel.barrenechea@strath.ac.uk}}.
}
}

\author{
  Alex Lukyanov
}
\address{
  Alex Lukyanov
  \thanks{
    Department of Mathematics and Statistics,
    University of Reading,
    {\tt{a.lukyanov@reading.ac.uk}}.
}}

\author{
  Tristan Pryer
}
\address{
  Tristan Pryer
  \thanks{
    Department of Mathematical Sciences,
    University of Bath, Bath BA2 7AY, UK
    {\tt{tmp38@bath.ac.uk}}.
}}

\author{
  Alex Trenam
}
\address{
  Alex Trenam 
  \thanks{
    The Maxwell Institute for Mathematical Sciences \& Mathematics-Driven Innovation Centre, School of Mathematical and Computer Sciences, Heriot-Watt University, Edinburgh, EH14 4AS, UK.    
    {\tt{A.Trenam@hw.ac.uk}}.
}}

\title{A finite element method for a non-Newtonian dilute polymer fluid}
\date{\today}

\begin{document}

\maketitle
\begin{abstract}
  We study the discretisation of a uniaxial (rank-one) reduction of
  the Oldroyd-B model for dilute polymer solutions, in which the
  conformation tensor is represented as $\sig = \vec b \otimes \vec
  b$. Building on structural analogies with MHD, we formulate a finite
  element framework compatible with the de Rham complex, so that the
  discrete velocity is exactly divergence-free. The spatial
  discretisation combines an interior-penalty treatment of viscosity
  with upwind transport to control stress layers and we prove inf-sup
  conditions on the mixed pairs. For time-stepping, we design an IMEX
  scheme that is linear at each step and show well-posedness of the
  fully discrete problem together with a discrete energy law mirroring
  the continuum dissipation. Numerical experiments on canonical
  benchmarks (lid-driven cavity, pipe-with-cavity and $4{:}1$ planar
  contraction) demonstrate accuracy and robustness for
  moderate-to-high Weissenberg numbers, capturing sharp stress
  gradients and corner singularities while retaining the efficiency
  gains of the uniaxial model. The results indicate that de
  Rham-compatible discretisations coupled with energy-stable IMEX time
  integration provide a reliable pathway for viscoelastic computations
  at elevated elasticity.
\end{abstract}

\section{Introduction}

Viscoelastic fluids, such as polymer solutions, exhibit complex
behaviours that differ significantly from those of Newtonian fluids
due to the presence of long-chain molecules whose microstructure
meaningfully affects macroscopic flow. Industrial examples include
molten plastics and lubricants, while blood is a canonical biological
example. Unlike Newtonian flows, governed solely by the Navier-Stokes
equations, viscoelastic flows require an additional constitutive
relation to account for polymer-induced stresses, leading to systems
of partial differential or integro-differential equations
\cite{hinch2021oldroyd,bird1987dynamics,owens2002computational}. Among
many constitutive models, the Oldroyd-B model is a widely used
benchmark that couples viscous and elastic responses by combining a
Newtonian solvent with a Maxwell-type polymer contribution
\cite{fattal2004constitutive,hinch2021oldroyd}. Despite its
simplicity, Oldroyd-B captures key aspects of viscoelastic behaviour
and remains central in computational rheology
\cite{fattal2005time,boyaval2010lid}. Flow regimes are commonly
organised by the Weissenberg and Deborah numbers, the Weissenberg
number $\Wi$ quantifies the significance of nonlinear normal-stress
effects \cite{dealy2010weissenberg,pakdel1997cavity}, whereas the
Deborah number $\De$ measures the dominance of elastic over viscous
timescales \cite{bird1987dynamics}. At low $\De$ one recovers nearly
Newtonian behaviour, whereas at high $\Wi$ purely elastic phenomena,
including elastic turbulence and flow-induced instabilities, can
appear and complicate both analysis and computation
\cite{pakdel1998cavity}.

A central numerical difficulty is the \emph{high Weissenberg number
problem} (HWNP). As $\Wi$ increases, steep stress gradients and
exponential boundary layers form, which frequently trigger loss of
convergence or spurious oscillations \cite{pan2009simulation}. The
Oldroyd-B model also features unbounded extensional viscosity, further
exacerbating instability \cite{hinch2021oldroyd}. Writing the
constitutive law
\begin{equation}\label{eq:constitutive_law}
  \frac{\partial \sig}{\partial t} + (\vec u \cdot \grad)\sig
  - (\grad \vec u)\sig - \sig(\grad \vec u)^{T}
  = - \frac{1}{\Wi} (\sig - \vec I)
\end{equation}
as an ODE along characteristics induced by the deformation gradient
clarifies the mechanism by which eigenvalues of $\sig$ can grow and
precipitate numerical breakdown \cite{ashby2025discretisation}. This
combination of physical stiffness and geometric transport underlies
the need for structure-preserving discretisations.

This work introduces a finite element framework that addresses the
computational challenges above by combining a model reduction with
compatible spatial discretisation and energy-stable time
integration. The starting point is a uniaxial (rank-one) reduction of
Oldroyd-B, originally proposed in \cite{fouxon2003spectra}, in which
the conformation tensor is approximated by a dyad, thereby retaining
dominant alignment and stretch at reduced cost in multiple
dimensions. Inspired by magnetohydrodynamics (MHD), we develop finite
element methods that respect the de~Rham complex: velocity and
pressure are approximated in $H(\mathrm{div};\Omega)$- and
$L^{2}$-conforming spaces to enforce incompressibility exactly at the
discrete level, while the polymeric vector field is placed in an
$H(\mathrm{curl};\Omega)$-conforming space to mirror the
transport-rotation structure. We establish well-posedness of the and
discrete energy laws mirroring the continuum dissipation for both the
semi- and fully-discrete schemes. Numerical experiments on demanding
benchmarks, including lid-driven cavity, pipe-with-cavity and planar
contraction flows, confirm robustness at high $\Wi$ and demonstrate
improved computational efficiency relative to traditional approaches.

The HWNP has been recognised since early computational rheology
studies \cite{crochet2012numerical,keunings1986high,boger1987viscoelastic,oliveira2007effect},
where steady solutions in benchmark problems (e.g. contractions and
expansions) often failed at modest $\Wi$. Subsequent developments
pursued two complementary directions, stabilising formulations and
structure-preserving discretisations. On the formulation side, the
log-conformation representation \cite{fattal2004constitutive,fattal2005time}
guarantees positive-definiteness of the conformation tensor and
dramatically extends attainable $\Wi$; square-root and
hyperbolic-tangent mappings provide related positivity-preserving
transforms \cite{lozinski2003energy,jafari2018property}, with
comparative assessments reported in \cite{zhang2021comparative}.
Classical stabilisations for advective transport, including
SUPG/Petrov-Galerkin ideas \cite{hughes1987new} and DEVSS-type
splittings \cite{guenette1995new}, have also been employed, as have
finite-volume and DG variants tailored to high P\'eclet transport of
stress \cite{comminal2015robust,cockburn2006introduction}. On the
structural side, energy and entropy estimates for Oldroyd-B and
related models \cite{lozinski2003energy,hu2007new,lions2000global,barrett2005existence}
have motivated schemes that ensure discrete free-energy decay
\cite{boyaval2010lid,boyaval2009free}. Our IMEX formulation inherits a
discrete energy law without resorting to fully nonlinear solves at
each time step, complementing implicit log-conformation approaches
\cite{hulsen2005flow}.

Compatibility of mixed spaces plays a central role in incompressible
and multiphysics flows. Finite element exterior calculus (FEEC)
provides a unifying perspective on stable discretisation through
commuting projections on the de~Rham complex \cite{arnold2006finite}.
For incompressible flow and MHD, $H(\mathrm{div})$-conforming
velocities and $H(\mathrm{curl})$-conforming fields improve stability
and preserve differential constraints discretely
\cite{boffi2013mixed,ogilvie2003relation}. Our use of
$H(\mathrm{div})/L^{2}$ for velocity-pressure eliminates spurious
pressure modes and provides exact enforcement of the divergence
constraint, while the $H(\mathrm{curl})$ discretisation for the
polymeric vector field emulates the curl-advection structure that
appears in both viscoelastic transport and MHD analogies
\cite{comminal2015robust}.  From the algorithmic viewpoint, IMEX time
integrators \cite{ascher1995implicit} are natural for stiff
transport-relaxation systems; here we show that, coupled with the
skew-symmetric treatment of convection and appropriate upwinding, they
yield a discrete energy decay consistent with the continuum law.

Reduced-order modelling aims to alleviate the high computational
burden of fully resolved viscoelastic simulations. Proper orthogonal
decomposition (POD) and Galerkin reduced models have been explored for
viscoelastic turbulence and complex flows \cite{chen2017pod}, often
combined with regularisation or stabilisation
\cite{chetry2023comparing}. At the model level, rank-limited (e.g.
rank-one) approximations of the conformation tensor provide a
physics-based reduction that captures dominant alignment and stretch
while preserving key couplings \cite{fouxon2003spectra}. Our
framework integrates such a uniaxial reduction with compatible
discretisation and energy stability, thereby achieving substantial
efficiency gains without sacrificing robustness.

Finally, the analogy between viscoelastic constitutive dynamics and
magnetic induction in MHD provides conceptual and practical guidance.
Ogilvie and Proctor \cite{ogilvie2003relation} formalised
correspondences between Oldroyd-B and MHD, including instability
analogues, which motivate the enforcement of structure-preserving
constraints (e.g. divergence-like and positivity conditions) and the
use of $H(\mathrm{div})/H(\mathrm{curl})$ spaces. Empirically,
viscoelastic instabilities in cavity and contraction flows
\cite{pakdel1997cavity,pakdel1998cavity,comminal2015robust,oliveira2007effect}
align with this viewpoint, and numerical treatments developed in MHD
inform our stabilisation and discretisation choices.

Section \ref{sec:prelim} presents the governing equations, notation
and basic assumptions. Section \ref{sec:modelreduction} introduces the
uniaxial reduction, discusses its modelling implications and
establishes key properties, including the relation to the continuum
energy law. Section \ref{sec:weak} develops the weak formulation,
identifies the bilinear and trilinear forms used in discretisation and
records identities needed for stability. Section \ref{sec:spatial}
details the finite element spaces, their de~Rham compatibility and
stabilised discrete operators, and proves the relevant inf-sup
conditions together with existence for the stationary problem. Section
\ref{sec:fulldis} presents the IMEX scheme and proves uniqueness per
time step together with the discrete energy law. Section
\ref{sec:numerics} reports numerical experiments for lid-driven
cavity, pipe-with-cavity and $4{:}1$ planar contraction flows. Section
\ref{sec:conclusion} summarises conclusions and outlines future
directions.

\section{Model Formulation and Preliminaries}
\label{sec:prelim}

In this section we fix notation, record the functional setting, and
state standing assumptions for the analysis to follow. Throughout,
$\Omega\subset\mathbb{R}^3$ denotes a bounded Lipschitz domain with
boundary $\partial\Omega$ and outward unit normal $\vec n$.

We write $\leb{p}(\omega)$ for the Lebesgue space on a measurable
subset $\omega\subset\mathbb{R}^d$, $d=1,2,3$, $1\le p\le\infty$, with norm
$\|\cdot\|_{\leb{p}(\omega)}$. For $p=2$ we use the inner product
$\langle\cdot,\cdot\rangle_\omega$ on $\leb{2}(\omega)$ and define
\begin{equation*}
  \leb{2}_0(\omega) := \{  q\in \leb{2}(\omega)  :  \langle q,1\rangle_\omega = 0  \}.
\end{equation*}
When $\omega=\Omega$ we suppress the subscript and write
$\langle\cdot,\cdot\rangle=\langle\cdot,\cdot\rangle_\Omega$.

For a nonnegative integer $k$, the Sobolev space $\sobh{k}(\omega)$ is defined as
\begin{equation}
  \sobh{k}(\omega)
  :=
  \{  \phi\in \leb{2}(\omega)  :  \D^{\vec\alpha}\phi \in \leb{2}(\omega)\ \text{for all}\ |\vec\alpha|\le k  \},
\end{equation}
where $\vec\alpha=(\alpha_1,\ldots,\alpha_d)$ is a multi-index,
$\D^{\vec\alpha}$ the associated weak derivative, and
$|\vec\alpha|=\sum_{i=1}^d\alpha_i$. We equip $\sobh{k}(\omega)$ with
the norm and seminorm
\begin{equation}
  \Norm{u}_{\sobh{k}(\omega)}^2 := \sum_{|\vec\alpha|\le k} \|\D^{\vec\alpha}u\|_{\leb{2}(\omega)}^2,
  \qquad
  \norm{u}_k^2 := \sum_{|\vec\alpha|=k} \|\D^{\vec\alpha}u\|_{\leb{2}(\omega)}^2.
\end{equation}
In particular
\begin{equation}
  \hoz(\Omega) := \{  \phi\in \sobh{1}(\Omega)  :  \phi|_{\partial\Omega}=0  \}.
\end{equation}

Vector- and tensor-valued fields are understood componentwise. For a
scalar field $u:\mathbb{R}^d\to\mathbb{R}$, $\grad u$ denotes the
gradient column vector; for a vector field $\vec
v:\mathbb{R}^d\to\mathbb{R}^d$, $\div\vec v$ is the scalar divergence
and $\grad\vec v$ the Jacobian matrix with entries $(\grad\vec
v)_{ij}=\D_j v_i$. The Laplacian is $\Delta u:=\div(\grad u)$, applied
componentwise to vectors. We employ the Frobenius inner product $A : B
:= \sum_{i,j}A_{ij}B_{ij}$ for matrices, and the associated norm
$\|A\|_{\mathrm F}^2:=A : A$. The curl of a vector field is denoted by
$\curl{\vec v}$.

For fields with divergence or curl in $\leb{2}$ we use the standard
Sobolev spaces
\begin{equation}
  \sobhdiv{\Omega}
  := \{  \vec\phi\in [\leb{2}(\Omega)]^d  :  \div\vec\phi \in \leb{2}(\Omega)  \}, \qquad
  \sobhcurl{\Omega}
  := \{  \vec\phi\in [\leb{2}(\Omega)]^d  :  \curl\vec\phi \in [\leb{2}(\Omega)]^d  \}.
\end{equation}
These spaces will host the discrete velocity (in $\sobhdiv{\Omega}$) and the
polymeric vector field (in $\sobhcurl{\Omega}$), respectively, so as
to reflect the underlying differential constraints.

Time-dependent problems are posed on $[0,T]$ with $T>0$. For a Banach
space $X$, the Bochner space $\leb{2}(0,T;X)$ is
\begin{equation}
  \leb{2}(0,T;X) := \big\{  u:[0,T]\to X \ \text{strongly measurable}  :  \int_0^T \|u(t)\|_X^2 \mathrm dt < \infty  \big\}.
\end{equation}
When needed we also employ $H^1(0,T;X^*)$ for time derivatives in a
dual space. We write $\|u\|_{\leb{2}(0,T;X)}^2 := \int_0^T
\|u(t)\|_X^2 \mathrm dt$ and, when $X=\sobh{k}(\Omega)$, abbreviate
$\leb{2}(0,T;\sobh{k})$.

Unless stated otherwise, all equalities and inequalities are to be
understood almost everywhere, and constants denoted by $C$ may change
from line to line but are independent of the mesh parameters
introduced later.

\section{Model Description and Uniaxial Reduction}
\label{sec:modelreduction}

In this section we formulate the governing equations for viscoelastic
fluid flow under the Oldroyd--B model and derive a reduced uniaxial
model. Let $\vec{u}:\Omega \times [0, T] \to \reals^d$ represent the
velocity field of an incompressible fluid in a domain $\Omega \subset
\reals^d$ with polyhedral boundary, and let $p:\Omega \times [0, T]
\to \reals$ denote the pressure. The viscoelastic effects are captured
by the symmetric \emph{conformation tensor} $\sig: \Omega \times [0,
T] \to \reals^{d \times d}$, describing the non-Newtonian contribution
to the stress.

The governing equations consist of the incompressible Navier-Stokes
equations coupled with a hyperbolic constitutive law for the evolution
of the conformation tensor. The Navier-Stokes equations for the
velocity and pressure are
\begin{equation}\label{eq:navier_stokes_modelreduction}
  \begin{split}
    \Re \left( \frac{\partial \vec{u}}{\partial t} + (\vec{u} \cdot \grad) \vec{u} \right)
    &=
    - \grad p + 2\beta \div\qp{\symgrad{\vec{u}}}
    + \frac{1 - \beta}{\Wi} \div \sig, \\
    \div \vec{u} &= 0, \\
    \vec u\vert_{\partial \Omega} &= \vec 0,
  \end{split}
\end{equation}
where the rate-of-strain tensor $\symgrad{\vec{u}}$ is defined as
\begin{equation}
  \symgrad{\vec{u}} = \frac{\grad \vec{u} + (\grad \vec{u})^T}{2}.
\end{equation}
The conformation tensor $\sig$ evolves according to the Oldroyd--B
constitutive law
\begin{equation}\label{eq:oldroyd_b_constitutive}
  \frac{\partial \sig}{\partial t} + (\vec{u} \cdot \grad) \sig
  - (\grad \vec{u}) \sig - \sig (\grad \vec{u})^T
  =
  -\frac{1}{\Wi} (\sig - \vec{I}),
\end{equation}
where $\vec{I}$ is the identity tensor.

The dimensionless parameters in the equations are the Reynolds number,
the Weissenberg number and the elasticity ratio. The Reynolds number
$\Re = \frac{\rho U L}{\eta}$ is the ratio of inertial to viscous
forces, where $\rho$ is the fluid density, $U$ a characteristic
velocity, $L$ a characteristic length and $\eta$ the dynamic
viscosity. The Weissenberg number $\Wi = \frac{\lambda U}{L}$
represents the ratio of the characteristic polymer relaxation time
$\lambda$ to the flow timescale $L / U$. The elasticity ratio
$\beta = \frac{\eta_s}{\eta}$ is the ratio of the solvent viscosity
$\eta_s$ to the total viscosity $\eta = \eta_s + \eta_p$, where
$\eta_p$ is the polymer contribution to the viscosity. These
parameters characterise the relative importance of inertial, viscous
and elastic effects in viscoelastic flows and govern the dynamic
behaviour of the system.

\begin{Rem}[Boundary and initial conditions]
  Throughout this work we assume no-slip boundary conditions for the
  velocity field $\vec{u}$, so that $\vec u = \vec 0$ on
  $\partial \Omega$. Under this assumption there is no inflow boundary
  for the conformation tensor $\sig$, and no explicit boundary
  conditions are required for $\sig$.

  In situations where the velocity boundary conditions are not
  trivial, for example when Dirichlet data prescribe a non-zero
  velocity $\vec{u} = \vec{g}$ on a portion
  $\Gamma \subset \partial \Omega$, suitable boundary conditions must
  also be imposed for the conformation tensor. The equation
  \eqref{eq:oldroyd_b_constitutive} is hyperbolic in space and
  therefore requires boundary data only on the inflow boundary
  \begin{equation}
    \Gamma_- = \{ \vec{x} \in \Gamma : \vec{g}(\vec{x}) \cdot \vec{n}(\vec x) < 0 \},
  \end{equation}
  where $\vec{n}$ is the outward unit normal to $\partial \Omega$. On
  $\Gamma_-$ we prescribe the inflow conformation via
  \begin{equation}
    \sig = \sig_{\text{in}} \quad \text{on } \Gamma_-,
  \end{equation}
  where $\sig_{\text{in}}$ is a positive-definite tensor representing
  the polymeric conformation entering the domain.

  At the initial time we set
  \begin{equation}
    \vec{u}(\vec x,0)=\vec{u}_0(\vec x), \qquad \sig(\vec x, 0) = \sig_0(\vec x),
  \end{equation}
  where $\vec u_0$ is the initial velocity (assumed solenoidal and
  compatible with the boundary conditions) and $\sig_0$ is a
  positive-definite tensor field representing the initial conformation
  state.
\end{Rem}

\subsection{Strong-stretch regime and uniaxial reduction}

As mentioned in the introduction, the analysis in
\cite{fouxon2003spectra} focuses on regimes with strongly stretched
polymers and large Weissenberg number, in which the conformation
tensor dominates the isotropic contribution, $\|\sig\|\gg 1$. In this
setting the relaxation term in \eqref{eq:oldroyd_b_constitutive} is
approximated by
\[
  -\Wi^{-1}(\sig - \vec I) \approx -\Wi^{-1}\sig,
\]
so that the constitutive law is replaced by the upper-convected
Maxwell-type approximation
\begin{equation}\label{eq:ucm_constitutive}
  \frac{\partial \sig}{\partial t} + (\vec{u} \cdot \grad) \sig
  - (\grad \vec{u}) \sig - \sig (\grad \vec{u})^T
  =
  -\frac{1}{\Wi} \sig.
\end{equation}
This step changes the isotropic equilibrium from $\sig = \vec I$ to
$\sig = \vec 0$, in the absence of flow the polymer stress relaxes to
zero, corresponding to an effectively Newtonian response.

On top of this strong-stretch approximation we adopt a uniaxial ansatz
for the conformation tensor,
\begin{equation}\label{eq:uniaxial_assumption}
  \sig = \vec{b} \otimes \vec{b},
\end{equation}
where $\vec{b}:\Omega \times [0, T] \to \reals^d$ is a vector field
representing the primary direction and magnitude of polymeric stress,
and $\otimes$ denotes the dyadic product. Physically, $\vec{b}$
encodes the dominant direction and magnitude of polymer chain
alignment and stretching within the fluid: its magnitude reflects the
extent of polymeric elongation, while its direction corresponds to the
principal axis along which polymer chains are aligned by the flow.

This assumption reflects a dominant alignment of polymer molecules in
the direction of $\vec{b}$, which is a reasonable approximation in
flows with strong uniaxial stretching or pronounced alignment of
stress. It is particularly relevant in regimes where polymer chains
predominantly align in a single direction due to flow-induced
stretching, such as in extensional or high-shear flows. One of our
goals is to explore, through computation, the range of applicability
of this modelling assumption.

Substituting \eqref{eq:uniaxial_assumption} into the constitutive
equation \eqref{eq:oldroyd_b_constitutive}, we rewrite the evolution
of $\sig$ in terms of $\vec{b}$. Differentiating $\sig$ with respect
to time gives
\begin{equation*}
  \frac{\partial \sig}{\partial t}
  =
  \frac{\partial \vec{b}}{\partial t} \otimes \vec{b}
  +
  \vec{b} \otimes \frac{\partial \vec{b}}{\partial t}.
\end{equation*}
Similarly,
\begin{align*}
  (\vec{u} \cdot \grad) \sig
  &=
  \left((\vec{u} \cdot \grad) \vec{b}\right) \otimes \vec{b}
  +
  \vec{b} \otimes \left((\vec{u} \cdot \grad) \vec{b}\right),\\
  (\grad \vec{u}) \sig &= \left((\grad \vec{u}) \vec{b}\right) \otimes \vec{b}, \\
  \sig (\grad \vec{u})^T &= \vec{b} \otimes \left((\grad \vec{u}) \vec{b}\right).
\end{align*}
Each term remains a dyadic product of $\vec{b}$ with itself, so we can
equate the coefficients of these products. This yields the evolution
equation for $\vec{b}$:
\begin{equation}
  \label{eq:b_evolution}
  \frac{\partial \vec{b}}{\partial t} + (\vec{u} \cdot \grad) \vec{b}
  - (\vec{b} \cdot \grad) \vec{u} = -\frac{1}{\Wi} \vec{b}.
\end{equation}
The uniaxial assumption ensures that the tensor evolution remains
aligned with $\vec{b}$: the terms in
\eqref{eq:oldroyd_b_constitutive} act only to modify the magnitude and
orientation of $\vec{b}$, rather than introducing additional stress
directions or degrees of freedom.

\begin{Rem}[Uniaxial boundary and initial conditions]
  As for the full system, no boundary conditions are required for
  $\vec b$ in the absence of inflow boundaries, while on inflow
  portions of the boundary $\vec b$ must be prescribed consistently
  with the imposed conformation. 

  The initial condition for $\vec{b}$ can be derived from the initial
  conformation tensor $\sig_0$ by an eigenvalue decomposition. Since
  $\sig_{\text{in}}$ and $\sig_0$ are positive definite in physical
  flows, they are not rank-one, so an approximation is required to
  ensure compatibility with the reduced model. A natural choice is to
  retain only the dominant stress mode, extracting the largest
  eigenvalue $\lambda_1$ and corresponding eigenvector $\vec{v}_1$ so
  that
  \begin{equation}
    \sig_0 \approx \lambda_1 \vec{v}_1 \otimes \vec{v}_1.
  \end{equation}
  The reduced variable is then defined as $\vec{b} = \sqrt{\lambda_1}
  \vec{v}_1$. This projection ensures that the reduced model captures
  the dominant stress direction while remaining consistent with the
  original conformation tensor.
\end{Rem}

\begin{Rem}[The ``do nothing'' configurations]
  \label{rem:do-nothing}
  In the full Oldroyd--B model, the identity tensor $\vec{I}$
  represents the ``do nothing'' configuration for the conformation
  tensor $\sig$. The relaxation term $-\Wi^{-1} (\sig - \vec{I})$
  ensures that $\sig$ relaxes toward this isotropic equilibrium state
  in the absence of external forces or flow. In eigenvalue terms,
  $\sig = \vec{I}$ corresponds to all eigenvalues being equal to $1$,
  reflecting an isotropic stress state with no preferred direction.

  Under the strong-stretch approximation \eqref{eq:ucm_constitutive}
  and the uniaxial assumption \eqref{eq:uniaxial_assumption}, the
  conformation tensor is expressed as $\sig = \vec{b} \otimes \vec{b}$,
  where $\vec{b}$ encodes the polymeric stress direction and
  magnitude. This imposes a rank-one structure on $\sig$, so $\sig$
  can have at most one nonzero eigenvalue, namely $\|\vec{b}\|^2$, and
  the relaxation acts towards $\vec b = 0$ rather than towards
  $\sig = \vec I$. Consequently, the isotropic equilibrium state
  $\sig = \vec{I}$ of Oldroyd--B cannot be represented within the
  reduced model, whose natural ``do nothing'' configuration is
  $\vec{b} = \vec{0}$ (hence $\sig = \vec 0$), corresponding to
  vanishing polymer stress and an effectively Newtonian response.
\end{Rem}

The uniaxial assumption reduces the dimensionality of the system. The
original tensorial constitutive law requires evolving $d(d+1)/2$
independent components of the symmetric tensor $\sig$, while the
vector $\vec{b}$ has only $d$ components. This simplification is
particularly advantageous in three-dimensional problems, where the
number of degrees of freedom is reduced from six (for $\sig$) to three
(for $\vec{b}$).

This reduction is not without a price, however: there is an inherent
limit to the model's ability to capture secondary stress components
orthogonal to $\vec{b}$. This trade-off needs to be carefully
considered in applications involving complex flow patterns or strong
interactions between multiple stress components.

A further structural feature of the reduced model is that the vector
field $\vec b$ becomes asymptotically solenoidal, as the following
result shows.

\begin{Pro}[Asymptotic solenoidal stress vector field]
  Let $\vec{b}$ solve \eqref{eq:b_evolution} and assume that the
  velocity field $\vec{u}$ satisfies $\div \vec{u} = 0$. Then
  \begin{equation}
    \lim_{t \to \infty} \div \vec{b}(t,\vec x) = 0 \qquad \text{for all } \vec x\in\Omega.
  \end{equation}
\end{Pro}

\begin{Proof}
  Define $y = \div \vec{b}$. Taking the divergence of both sides of
  \eqref{eq:b_evolution}, and using $\div \vec{u} = 0$, we obtain
  \begin{equation}
    \frac{\partial y}{\partial t} + (\vec{u} \cdot \grad) y = -\frac{1}{\Wi} y.
  \end{equation}
  This is a linear first-order partial differential equation for $y$,
  in which the convective term $(\vec{u} \cdot \grad) y$ advects the
  divergence along fluid trajectories, and the relaxation term
  $-\frac{1}{\Wi} y$ drives $y$ toward zero.

  Along a characteristic curve $\vec{x}(t)$ of the velocity field,
  satisfying $\dot{\vec x}(t) = \vec u(\vec x(t),t)$, the quantity
  $y(t,\vec x(t))$ solves the ordinary differential equation
  \begin{equation}
    \frac{\mathrm d}{\mathrm dt} y(t,\vec x(t))
    = -\frac{1}{\Wi} y(t,\vec x(t)),
  \end{equation}
  with solution
  \begin{equation}
    y(t, \vec{x}(t)) = y(0, \vec{x}(0)) \exp \left(-\frac{t}{\Wi}\right).
  \end{equation}
  As $t \to \infty$, the exponential decay ensures that
  $y(t,\vec x(t)) \to 0$, and therefore $\div \vec{b}(t,\vec x) \to 0$
  along every trajectory, as claimed.
\end{Proof}

\begin{Rem}[Relaxation in the reduced model]
  The proposition shows that, under the combined effects of relaxation
  and advection, the vector field $\vec{b}$ asymptotically approaches
  a divergence-free configuration. In the absence of external forcing,
  the polymeric stress field therefore tends toward an equilibrium
  state in which stress magnitudes decay over time and no net
  divergence is generated. This is consistent with the expectation
  that, in the long-time limit, the reduced model relaxes toward a
  stable configuration governed by the balance between advection and
  relaxation.
\end{Rem}

Collecting the above simplifications and assumptions, and in
particular using the rank-one ansatz $\sig = \vec b\otimes\vec b$ and
the divergence-free constraints $\div \vec u = 0$ and
$\div \vec b = 0$ (motivated by the discussion and Proposition above),
the reduced system analysed in this work reads
\begin{equation}\label{eq:reduced_model_divfree}
  \begin{aligned}
    \Re \left( \frac{\partial \vec{u}}{\partial t}
      + (\vec{u} \cdot \grad) \vec{u} \right)
    - 2\beta \div\qp{\symgrad{\vec{u}}}
    + \grad p
    &=
     \frac{1 - \beta}{\Wi} (\vec{b} \cdot \grad) \vec{b}
     &&\text{in } \Omega \times (0,T), \\
    \frac{\partial \vec{b}}{\partial t}
    + (\vec{u} \cdot \grad) \vec{b}
    - (\vec{b} \cdot \grad) \vec{u}
    &=
    -\frac{1}{\Wi} \vec{b}
    &&\text{in } \Omega \times (0,T), \\
    \div \vec{u} &= 0
    &&\text{in } \Omega \times (0,T), \\
    \div \vec{b} &= 0
    &&\text{in } \Omega \times (0,T), \\
    \vec u &=\vec 0
    &&\text{on } \partial\Omega \times (0,T),
    \\
    \vec u(\vec x, 0) &=\vec u_0(\vec x)
    &&\text{in } \Omega \times \{0 \},
    \\
    \vec b(\vec x, 0) &=\vec b_0(\vec x)
    &&\text{in } \Omega \times \{0 \},
  \end{aligned}
\end{equation}

Note that the upper-convected derivative of $\vec{b}$ admits the
representation
\begin{equation}
  \label{eq:mhdrelation}
  \frac{\partial \vec{b}}{\partial t}
  +
  (\vec{u} \cdot \nabla) \vec{b}
  -
  (\vec{b} \cdot \nabla) \vec{u}
  =
  \frac{\partial \vec{b}}{\partial t}
  -
  \curl (\vec{u} \times \vec{b}),
\end{equation}
where we have used the vector identity
\[
  \curl(\vec u\times\vec b)
  = (\vec b\cdot\nabla)\vec u - (\vec u\cdot\nabla)\vec b
    + \vec u \div\vec b - \vec b \div\vec u
\]
together with $\div\vec u = \div\vec b = 0$. In particular, the
divergence-free constraints eliminate terms associated with isotropic
expansion or compression of $\vec b$, so that only advection and
shearing/stretching contribute to the evolution of $\vec b$.

The divergence-free property $\div \vec{b} = 0$ is consistent with the
interpretation of $\vec{b}$ as the dominant direction and magnitude of
polymeric stress, it indicates that there are no artificial sources or
sinks of stress within the domain, and that stress is transported and
redistributed by the flow without spurious accumulation.

This structure highlights a strong analogy with magnetohydrodynamics
(MHD), where the magnetic field is divergence-free and its evolution
contains a term $\curl (\vec{u} \times \vec{B})$ describing the
advection and stretching of magnetic field lines by the velocity
field. In the present viscoelastic setting, the term
$\curl (\vec{u} \times \vec{b})$ plays an analogous role, capturing
the advection and alignment of polymeric stress vectors. This analogy
provides useful insight into the dynamics of $\vec{b}$ and motivates
the adaptation of numerical techniques from MHD to viscoelastic fluid
simulations.

\section{Weak Formulation}
\label{sec:weak}

In this section we derive the weak formulation of the system described
in \S\ref{sec:modelreduction}. The weak formulation provides a
framework suitable for numerical approximation by finite element
methods.

To express the governing equations in weak form, we introduce the
following bilinear and trilinear forms:
\begin{equation}
  \begin{aligned}
    \mathcal{O}(\vec{w}; \vec{u}, \vec{v})
    &:= \left\langle (\vec{w} \cdot \nabla) \vec{u}  , \vec{v} \right\rangle, \\
    \mathcal{A}(\vec{u}, \vec{v})
    &:= \langle \nabla \vec{u} , \nabla \vec{v}\rangle, \\
    \mathcal{C}(\vec{v}, q)
    &:= -\langle  \div{\vec{v}}, q \rangle, \\
    \mathcal{D}(\vec{c}, s)
    &:= \langle \vec{c} , \nabla s \rangle.
  \end{aligned}
\end{equation}
Here $\mathcal{O}$ is the trilinear convective form and $\mathcal{A}$
is the bilinear form associated with viscous diffusion. The form
$\mathcal{C}$ couples velocity and pressure and is used to enforce the
divergence-free condition for the velocity field, while
$\mathcal{D}$ couples the polymeric vector field to a scalar
multiplier and will be used to impose the divergence constraint on
$\vec{b}$.

The solenoidal constraint $\div \vec b = 0$ motivates the introduction
of an additional Lagrange multiplier $r$ in order to close the
variational system and enforce this condition weakly. Since $r$ enters
the formulation only through its gradient, it is determined up to an
additive constant and we work in a quotient (or zero-mean) space. We
therefore take as function spaces for the velocity, pressure,
polymeric vector field and scalar multiplier
\begin{equation}
  \mathcal{U} := \sobh{1}_0(\Omega)^d,
  \quad
  \mathcal{P} := L^2_0(\Omega),
  \quad
  \mathcal{B} := \sobh{1}(\Omega)^d,
  \quad
  \mathcal{R} := \sobh{1}(\Omega) / \mathbb R
\end{equation}

\subsection{Variational formulation}

The weak formulation of the system is: for almost every
$t \in (0,T)$, seek $\left( \vec{u}, p, \vec{b}, r \right) \in
\mathcal{U} \times \mathcal{P} \times \mathcal{B} \times \mathcal{R}$
such that
\begin{equation}
  \label{eq:weakform}
  \begin{split}
    \Re \qp{\ltwop{\partial_t \vec{u}}{\vec{v}}
      +
      \mathcal{O}(\vec{u}; \vec{u}, \vec{v})}
    +
    \mathcal{C}(\vec{v}, p)
    &=
    -\beta \mathcal{A}(\vec{u}, \vec{v})
    \\
    &\quad +
    \frac{1 - \beta}{\Wi} \left( \frac{1}{2} \mathcal{O}(\vec{b}; \vec{b}, \vec{v})
      - \frac{1}{2} \mathcal{O}(\vec{b}; \vec{v}, \vec{b}) \right),
    \\
    \ltwop{\partial_t \vec{b}}{\vec{c}}
    +
    \mathcal{O}(\vec{u}; \vec{b}, \vec{c})
    -
    \frac{1}{2} \mathcal{O}(\vec{b}; \vec{u}, \vec{c})
    +
    \frac{1}{2} \mathcal{O}(\vec{b}; \vec{c}, \vec{u})
    +
    \mathcal{D}(\vec{c}, r)
    &=
    -\frac{1}{\Wi} \ltwop{\vec{b}}{\vec{c}},
    \\
    \mathcal{C}(\vec{u}, q) &= 0, \qquad \mathcal{D}(\vec{b}, s) = 0
  \end{split}
\end{equation}
for all test functions
$\left( \vec{v}, q, \vec{c}, s \right) \in \mathcal{U} \times
\mathcal{P} \times \mathcal{B} \times \mathcal{R}$.

\begin{Rem}[Symmetric gradient equivalence]
  In deriving the weak formulation, and in particular the viscous
  terms, we use the identity
  \begin{equation}
    \int_\Omega \nabla \vec{u} : \nabla \vec{v}  \mathrm d \vec x
    =
    \int_\Omega \vec \varepsilon(\vec{u}) : \vec \varepsilon(\vec{v})  \mathrm d \vec x
    + \frac{1}{2} \int_\Omega (\div \vec{u})(\div \vec{v})  \mathrm d \vec x,
  \end{equation}
  which holds for sufficiently smooth vector fields and appropriate
  boundary conditions. In particular, for divergence-free test or
  trial functions the second term vanishes and we obtain
  \begin{equation}
    \int_{\Omega} \vec \varepsilon(\vec{u}) : \vec \varepsilon(\vec{v})  \mathrm d \vec x
    =
    \int_{\Omega} \nabla \vec{u} : \nabla \vec{v}  \mathrm d \vec x
    =
    \mathcal{A}(\vec{u}, \vec{v}).
  \end{equation}
  This justifies our choice of $\mathcal{A}$ in terms of the full
  gradient rather than the symmetric gradient.
\end{Rem}

\begin{Rem}[Skew symmetrisation]
  At the continuous level, thanks to the fact that $\vec u$ is
  solenoidal and vanishes on the boundary, the convective form
  satisfies the antisymmetries
  \begin{equation}\label{skew-symmetry-Cont}
    \mathcal{O}(\vec{u}; \vec{u}, \vec{u})
    =
    \mathcal{O}(\vec{u}; \vec{b}, \vec{b}) = 0.
  \end{equation}
  These cancellations will be instrumental in the stability analysis
  below.

  At the discrete level, the choices made for the forms
  $\mathcal{C}$ and $\mathcal{D}$ ensure that the discrete velocity is
  exactly divergence-free, whereas the discrete approximation to
  $\vec b$ is not. To mitigate the possible negative impact of this
  asymmetry on stability, we recall that for any solenoidal vector
  field $\vec{w}$, and smooth vector fields $\vec u$ and $\vec v$
  (with at least one of them vanishing on the boundary), integration
  by parts yields
  \[
    \mathcal{O}(\vec{w}; \vec{u}, \vec{v})
    = - \mathcal{O}(\vec{w}; \vec{v}, \vec{u}),
  \]
  and hence
  \begin{equation}
    \label{eq:symmetry}
    \mathcal{O}(\vec{w}; \vec{u}, \vec{v})
    =
    \frac{1}{2} \Big(
      \mathcal{O}(\vec{w}; \vec{u}, \vec{v})
      -
      \mathcal{O}(\vec{w}; \vec{v}, \vec{u})
    \Big).
  \end{equation}
  This skew-symmetrised representation underlies the convective
  combinations used in our discretisation and preserves the key
  cancellation properties at the continuous level.
\end{Rem}

The weak form \eqref{eq:weakform} inherits the dissipative structure
of the underlying continuum model, as the following result shows.

\begin{The}[Energy estimate]
  Assume that $\vec{u}_0 \in L^2(\Omega)^d$ and $\vec{b}_0 \in
  L^2(\Omega)^d$ are the initial conditions for the velocity and
  polymeric vector field, respectively. Then any solution of the weak
  formulation \eqref{eq:weakform} satisfies the energy identity
  \begin{equation}
    \label{eq:energy}
    \ddt\qp{
      \frac{\Re}{2} \Norm{\vec{u}}^2_{\leb{2}(\Omega)}
      +
      \frac{1 - \beta}{2\Wi}\Norm{\vec{b}}^2_{\leb{2}(\Omega)}
    }
    =
    -\beta \Norm{\nabla \vec{u}}^2_{\leb{2}(\Omega)}
    -\frac{1-\beta}{\Wi^2} \Norm{\vec{b}}^2_{\leb{2}(\Omega)}.
  \end{equation}
\end{The}

\begin{Proof}
  To derive \eqref{eq:energy}, we take $\vec{v} = \vec{u}$ in the
  first equation of \eqref{eq:weakform} and
  $\vec{c} = \frac{1-\beta}{\Wi} \vec{b}$ in the second one. With
  $q = p$ and $s = r$, adding the resulting equations gives
  \begin{multline}
    \Re \qp{\ltwop{\partial_t \vec{u}}{\vec{u}}
      +
      \mathcal O(\vec u;\vec u,\vec u)}
    +
    \frac{1 - \beta}{\Wi}\ltwop{\partial_t \vec{b}}{\vec{b}}
    +
    \frac{1-\beta}{\Wi}
    \mathcal O(\vec u;\vec b,\vec b)
    -
    \frac{1-\beta}{\Wi}
    \qp{
      \frac 12 \mathcal O(\vec b;\vec u,\vec b)
      -
      \frac 12 \mathcal O(\vec b;\vec b,\vec u)
    }
    \\
    =
    -
    \beta \mathcal A(\vec u, \vec u)
    +
    \frac{1-\beta}{\Wi}
    \qp{
      \frac 12 \mathcal O(\vec b;\vec b,\vec u)
      -
      \frac 12 \mathcal O(\vec b;\vec u,\vec b)
    }
    -
    \frac{1-\beta}{\Wi^2} \ltwop{\vec b}{\vec b}.
  \end{multline}
  The convective terms with first argument $\vec u$ vanish by the
  skew-symmetry \eqref{skew-symmetry-Cont}, since $\vec{u}$ is
  divergence-free and satisfies the homogeneous boundary conditions.
  The remaining convective terms, with first argument $\vec b$,
  cancel pairwise between the right- and left-hand sides. Using the
  definitions of $\mathcal A$ and the inner product then yields
  \eqref{eq:energy}, which completes the proof.
\end{Proof}

\renewcommand{\U}{\cU}
\renewcommand{\P}{\cP}
\renewcommand{\B}{\cB}
\renewcommand{\R}{\cR}

\section{Spatial Discretisation}
\label{sec:spatial}

In this section we introduce the spatial discretisation of the
problem, using finite element spaces chosen to be compatible with the
energy argument above and with the de~Rham complex. Let
$\mathcal{T}_h$ be a regular subdivision of the domain $\Omega$ into
disjoint simplicial elements $K$ (triangles in $d=2$, tetrahedra in
$d=3$). By $\mathcal{E}_h$ we denote the set of all interior
$(d-1)$-dimensional facets (edges in 2D, faces in 3D) associated with
the subdivision $\mathcal{T}_h$. We assume that the subdivision is
shape-regular, meaning that the ratio of the diameters of neighbouring
elements is uniformly bounded.

For each element $K \in \mathcal{T}_h$ let $\mathcal{P}_p(K)$ denote
the polynomials on $K$ of total degree at most $p$, and define the
broken polynomial space
\begin{equation}
  \mathcal{P}_p(\mathcal{T}_h)
  :=
  \left\{ u \in L^2(\Omega) : u|_K \in \mathcal{P}_p(K) \ \forall K \in \mathcal{T}_h \right\}.
\end{equation}
We write $h_K := \mathrm{diam}(K)$ for the element diameter and
collect these into the elementwise constant function
$h:\Omega \to \mathbb{R}$, with $h|_K = h_K$. For an interior facet
$e = K^+\cap K^- \subset \mathcal{E}_h$ we set
$h|_e = (h_{K^+} + h_{K^-})/2$, while for a boundary facet
$e \subset \partial K \cap \partial \Omega$ we set $h|_e = h_K$.

For a vector field $\vec{\phi}$, let $K^+$ and $K^-$ denote elements
sharing an interior facet $e \subset \mathcal{E}_h$, with outward unit
normals $\vec{n}^+$ and $\vec{n}^-$. The tensor jump and average
operators are defined, respectively, as
\begin{align}
  &\tjump{\vec{\phi}} := \vec{\phi}^+\otimes\vec{n}^+ + \vec{\phi}^-\otimes\vec{n}^-,
  &&\avg{\vec{\phi}} := \half\qp{\vec{\phi}^+ + \vec{\phi}^-},
\end{align}
where $\vec{\phi}^\pm$ are the traces of $\vec{\phi}$ on $e$ from
$K^\pm$. For boundary facets $e \subset \partial \Omega\cap K^+$ we
set
\begin{align}
  &\tjump{\vec{\phi}} := \vec{\phi}^+,
  &&\avg{\vec{\phi}} := \vec{\phi}^+.
\end{align}
Given an advecting field (typically the velocity), we choose the
labelling so that $K^+$ is the upwind element and $K^-$ the downwind
one on $e$. The upwind jump of $\vec{\phi}$ across $e$ is then
\begin{equation}
  \ujump{\vec \phi} := \vec\phi^+ - \vec\phi^-.
\end{equation}

\subsection{Finite element spaces}

We use finite element spaces compatible with the de~Rham complex to
discretise the vector fields, scalar fields and operators. In view of
\eqref{eq:mhdrelation}, for the polymeric vector field $\vec{b}$ we
adopt a curl-conforming discretisation based on the N\'ed\'elec
finite element spaces of the second kind:
\begin{equation}
  \mathcal{B}_h
  :=
  \mathbb{NED}_p(\Omega, \mathcal{T}_h)
  =
  \left\{ \vec{c} \in H(\mathrm{curl}; \Omega) : \vec{c}|_K \in \mathcal{P}_p(K)^d,\ \forall K \in \mathcal{T}_h \right\}.
\end{equation}
Approximating $\vec b$ in $\mathcal{B}_h$ yields a nonconforming
approximation with respect to $H^1(\Omega)^d$, since functions in
$\mathcal{B}_h$ have continuous tangential components but may be
discontinuous in their normal components across element interfaces.

For the velocity field $\vec{u}$ we employ a divergence-conforming
discretisation using the Brezzi-Douglas-Marini finite element spaces:
\begin{equation}
  \mathcal{U}_h :=
  \mathbb{BDM}_p(\Omega, \mathcal{T}_h) =
  \left\{
    \vec{v} \in H(\mathrm{div}; \Omega) :
     \vec{v}|_K \in \mathcal{P}_p(K)^d \ \text{and} \ \div \vec{v} \in \mathcal{P}_{p-1}(K), \ \forall K \in \mathcal{T}_h
  \right\}.
\end{equation}
The velocity is approximated in $\mathcal{U}_h$, so that the normal
component of the discrete velocity is continuous across facets, while
the tangential component may be discontinuous. Homogeneous velocity
boundary conditions are not imposed strongly in $\mathcal{U}_h$, they
will instead be enforced weakly by Nitsche-type terms in the
variational formulation below.

For the pressure $p$ and scalar field $r$, we use standard
discontinuous and continuous polynomial spaces, with a zero-mean
constraint on the latter:
\begin{equation}
  \begin{split}
    \mathcal{P}_h &:= \left\{ q_h \in L^2_0(\Omega) : q_h|_K \in \mathcal{P}_{p-1}(K), \ \forall K \in \mathcal{T}_h \right\}, \\
    \mathcal{R}_h &:= \left\{ s_h \in H^1(\Omega) : s_h|_K \in \mathcal{P}_p(K) \ \forall K \in \mathcal{T}_h,\ \langle s_h,1\rangle_\Omega = 0 \right\}.
  \end{split}
\end{equation}
On $\mathcal{R}_h$ the seminorm $\|s_h\|_{\mathcal R_h} := \|\nabla s_h\|_{L^2(\Omega)}$
is therefore a norm, by Poincar\'e-Wirtinger.

\subsection{Discrete de Rham complex}

The chosen finite element spaces fit into a discrete de~Rham complex,
equipped with projection operators that commute with the differential
operators:
\begin{equation}
  \begin{tikzcd}[row sep=large, column sep=huge]
    H^1(\Omega) \ar[r, "\nabla"] \ar[d, "\Pi_{\mathcal{R}_h}" left]
    & H(\mathrm{curl}; \Omega) \ar[r, "\mathrm{curl}"] \ar[d, "\Pi_{\mathcal{B}_h}" left]
    & H(\mathrm{div}; \Omega) \ar[r, "\mathrm{div}"] \ar[d, "\Pi_{\mathcal{U}_h}" left]
    & L^2(\Omega) \ar[d, "\Pi_{\mathcal{P}_h}" left] \\
    \mathcal{R}_h \ar[r, "\nabla"]
    & \mathcal{B}_h \ar[r, "\mathrm{curl}"]
    & \mathcal{U}_h \ar[r, "\mathrm{div}"]
    & \mathcal{P}_h
  \end{tikzcd}
\end{equation}
Here $\Pi_{\mathcal{R}_h}$, $\Pi_{\mathcal{B}_h}$, $\Pi_{\mathcal{U}_h}$
and $\Pi_{\mathcal{P}_h}$ denote suitable bounded projection operators
such that, for instance,
$\mathrm{curl} \Pi_{\mathcal{B}_h} = \Pi_{\mathcal{U}_h} \mathrm{curl}$
on $H(\mathrm{curl};\Omega)$ and
$\mathrm{div} \Pi_{\mathcal{U}_h} = \Pi_{\mathcal{P}_h} \mathrm{div}$
on $H(\mathrm{div};\Omega)$. This structure ensures that the discrete
spaces preserve the compatibility properties of the continuous
operators.

\begin{Rem}[Nonconforming vector spaces]
  While the discrete spaces for the pressure and Lagrange multiplier,
  $\mathcal{P}_h \subset \mathcal{P}$ and $\mathcal{R}_h \subset
  \mathcal{R}$, are conforming to their respective continuous spaces,
  the vector-valued spaces $\mathcal{U}_h$ and $\mathcal{B}_h$ are not
  conforming with respect to the $H^1$-based spaces used in the weak
  formulation. Specifically,
  \[
    \mathcal{U}_h \subset H(\mathrm{div};\Omega)
    \quad\text{whereas}\quad
    \mathcal{U} = \sobh{1}_0(\Omega)^d,
    \qquad
    \mathcal{B}_h \subset H(\mathrm{curl};\Omega)
    \quad\text{whereas}\quad
    \mathcal{B} = \sobh{1}(\Omega)^d.
  \]
  This choice is deliberate, the pair
  $\mathcal{U}_h\times \mathcal{P}_h$ is inf-sup stable (see below),
  and, moreover, the de~Rham structure implies
  $\div \mathcal{U}_h = \mathcal{P}_h$, so that the discrete velocity
  is exactly divergence-free. The nonconforming but de~Rham-compatible
  choice of $\mathcal{U}_h$ and $\mathcal{B}_h$ is crucial in
  preserving the structural properties of the continuous problem and
  plays an important role in the stability analysis of the numerical
  method.
\end{Rem}

\subsection{Stabilised discrete forms}

To account for the nonconforming nature of the velocity space, we
introduce a stabilised bilinear form for the diffusive term in the
weak formulation.

In addition to the interior facets $\mathcal{E}_h$ introduced above,
let $\mathcal{E}_h^\partial$ denote the set of boundary facets
$e \subset \partial\Omega$. We write
$\mathcal{E}_h^\ast := \mathcal{E}_h \cup \mathcal{E}_h^\partial$ for
the collection of all facets, interior and boundary. With the jump and
average operators extended to boundary facets as in the previous
subsection, we define, for $\vec{u}_h,\vec{v}_h \in \mathcal U_h$,
\begin{equation}
  \mathcal{A}_h(\vec{u}_h, \vec{v}_h)
  :=
  \sum_{K \in \mathcal{T}_h}
    \int_K \nabla \vec{u}_h : \nabla \vec{v}_h  \,\mathrm d\vec x
    +
    \sum_{e\in\mathcal{E}_h^\ast}\int_{e}
      \Big(
        \frac{\sigma_e}{h} \tjump{\vec{u}_h} : \tjump{\vec{v}_h}
        - \avg{\nabla \vec{u}_h} : \tjump{\vec{v}_h}
        - \avg{\nabla \vec{v}_h} : \tjump{\vec{u}_h}
      \Big) \,\mathrm ds,
\end{equation}
where $\sigma_e > 0$ is a facet-wise stabilisation parameter
penalising tangential discontinuities. On boundary facets the
definition of the jump and average reduces to
$\tjump{\vec u_h} = \vec u_h\otimes\vec n$ and
$\avg{\nabla\vec u_h} = \nabla\vec u_h$, so that the same form
provides a Nitsche-type weak enforcement of the homogeneous velocity
boundary condition. In practice we take
\[
  \sigma_e =
  \begin{cases}
    \sigma_{\mathrm{i}}, & e \text{ interior},\\[0.3em]
    \sigma_{\mathrm{b}}, & e \subset \partial\Omega,
  \end{cases}
\]
with constants $\sigma_{\mathrm{i}},\sigma_{\mathrm{b}} > 0$.

\begin{Lem}[DG norm and coercivity]
  The DG norm $\|\cdot\|_{1,h}$ for a discrete velocity field
  $\vec{v}_h \in \mathcal{U}_h$ is defined as
  \begin{equation}
    \Norm{\vec{v}_h}_{1,h}^2
    :=
    \sum_{K\in\mathcal T_h}\Norm{\nabla \vec{v}_h}_{L^2(K)}^2
    +
    \sum_{e \in \mathcal{E}_h^\ast} \frac{1}{h} \Norm{\tjump{\vec{v}_h}}_{L^2(e)}^2.
  \end{equation}
  The bilinear form $\mathcal{A}_h$ is coercive with respect to this
  norm: for $\sigma$ chosen large enough, and for any
  $\vec{u}_h \in \mathcal{U}_h$, it holds that (see, for example,
  \cite{arnold2002unified})
  \begin{equation}
    \mathcal{A}_h(\vec{u}_h, \vec{u}_h) \geq \alpha \|\vec{u}_h\|_{1,h}^2,
  \end{equation}
  where $\alpha > 0$ is a constant independent of $h$.
\end{Lem}

Below we state the inf-sup conditions associated with the pairs
$\mathcal{U}_h\times\mathcal{P}_h$ and $\mathcal{B}_h\times\mathcal{R}_h$.
Although their proofs are standard, and follow from the results in
\cite[Section~2.5]{boffi2013mixed} (see also
\cite[Lemma~7 and Proposition~10]{HL02}), we summarise them here for
completeness.

\begin{Lem}[Inf-sup conditions]
  \label{lem:infsup}
  Let $(\mathcal{U}_h, \mathcal{P}_h)$ and $(\mathcal{B}_h,
  \mathcal{R}_h)$ denote the finite element spaces for the
  velocity-pressure and polymeric vector field-Lagrange multiplier
  pairs, respectively. Then:
  \begin{enumerate}
  \item There exists a constant $C_1 > 0$, independent of $h$, such
    that    
    \begin{equation}\label{Eq:inf-sup-C}
      \sup_{\vec{v}_h \in \mathcal{U}_h} \frac{\mathcal{C}(\vec{v}_h, q_h)}{\Norm{\vec{v}_h}_{1,h}}
      \geq
      C_1 \Norm{q_h}_{\leb2(\Omega)} \Foreach q_h \in \mathcal{P}_h.
    \end{equation}
    
  \item For all $s_h \in \mathcal{R}_h$ it holds that
    \begin{equation}
      \sup_{\vec{c}_h \in \mathcal{B}_h} \frac{\mathcal{D}(\vec{c}_h, s_h)}{\Norm{\vec{c}_h}_{\leb2(\Omega)}}
        \geq
        \Norm{\nabla s_h}_{\leb 2(\Omega)}.
    \end{equation}
  \end{enumerate}
  In particular, since $\mathcal R_h$ consists of zero-mean functions,
  the quantity $\Norm{\nabla s_h}_{\leb 2(\Omega)}$ defines a norm on
  $\mathcal R_h$ by the Poincar\'e-Wirtinger inequality.
\end{Lem}

\begin{Proof}
  For the first inf-sup condition, let $q_h \in \mathcal{P}_h$. By the
  standard surjectivity of the divergence operator on $L^2_0(\Omega)$,
  there exists $\vec{v} \in \sobh1_0(\Omega)^d$ such that
  $\div\vec v=q_h$ in $\Omega$ and
  $\|\vec v\|_{H^1(\Omega)}\le C \|q_h\|_{L^2(\Omega)}$ for some
  constant $C>0$ independent of $h$. Let
  $\Pi_{\mathcal U_h}: \sobh1_0(\Omega)^d \to \mathcal{U}_h$ denote a
  Fortin operator satisfying the commuting diagram property
  \begin{equation}
    \div (\Pi_{\mathcal U_h} \vec{v}) = \Pi_{\mathcal P_h} (\div \vec{v}),
  \end{equation}
  where $\Pi_{\mathcal P_h}$ is the $L^2$-projection onto $\mathcal
  P_h$. Since $q_h\in\mathcal P_h$ and $\mathcal P_h$ consists of
  functions with zero mean, $\Pi_{\mathcal P_h} q_h = q_h$, and thus
  \begin{equation}
    \mathcal{C}(\Pi_{\mathcal U_h} \vec{v}, q_h)
    =
    \mathcal{C}(\vec{v}, q_h)
    =
    \|q_h\|_{L^2(\Omega)}^2.
  \end{equation}
  Moreover, by the $H^1(K)$-stability of $\Pi_{\mathcal U_h}$ (see
  \cite[Proposition~2.5.1]{boffi2013mixed}) we obtain
  \begin{equation}
    \Norm{\Pi_{\mathcal U_h} \vec{v}}_{1,h}
    \leq
    C \Norm{\vec{v}}_{\sobh1(\Omega)}
    \leq C \Norm{q_h}_{\leb2(\Omega)},
  \end{equation}
  with a constant $C$ independent of $h$. Taking
  $\vec v_h = \Pi_{\mathcal U_h}\vec v$ in \eqref{Eq:inf-sup-C} yields
  the desired bound for $\mathcal C$.

  For the second inf-sup condition, let $s_h \in \mathcal{R}_h$ be
  given and define
  \begin{equation}
    \vec{c}_h := \nabla s_h.
  \end{equation}
  By the discrete de~Rham complex we have
  $\vec c_h \in \mathcal B_h$, and moreover
  \begin{equation}
    \Norm{\vec{c}_h}_{\leb2(\Omega)} = \Norm{\nabla s_h}_{\leb2(\Omega)},
    \qquad
    \mathcal{D}(\vec{c}_h, s_h)
    =
    \Norm{\nabla s_h}_{\leb2(\Omega)}^2.
  \end{equation}
  Hence
  \begin{equation}
    \frac{\mathcal{D}(\vec{c}_h, s_h)}{\Norm{\vec{c}_h}_{\leb2(\Omega)}}
    =
    \Norm{\nabla s_h}_{\leb2(\Omega)}.
  \end{equation}
  Taking the supremum over all $\vec{c}_h \in \mathcal{B}_h$ gives the
  asserted bound for $\mathcal D$, and the proof is complete.
\end{Proof}

For $\vec{w}_h$, $\vec{u}_h$ and $\vec{v}_h$, the upwind
discretisation of the convective term is defined as
\begin{equation}
  \mathcal{O}_h\qp{\vec{w}_h; \vec{u}_h, \vec{v}_h}
  :=
  \sum_{K \in \mathcal{T}_h}
    \int_K\qp{\vec{w}_h\cdot\nabla}\vec{u}_h\cdot\vec{v}_h \mathrm d \vec{x}
  -\int_{\mathcal{E}_h}\qp{
    \qp{\vec{w}_h\cdot\normal}\avg{\vec{u}_h}\cdot\ujump{\vec{v}_h}
    - \frac{\mu}{2}\lvert\vec{w}_h\cdot\normal\rvert\ujump{\vec{u}_h}\cdot\ujump{\vec{v}_h}
  }\mathrm ds,
\end{equation}
where $\vec{n}$ is the unit normal on each interior facet and
$\mu>0$ is an upwind parameter (in the computations below we take
$\mu=1$).

We now summarise stability and continuity properties of the upwind
convective form. For proofs we refer to, for example,
\cite[Propositions~4.2 and 4.3]{CKS05}.

\begin{The}[Properties of $\mathcal{O}_h$]\label{Theo:Properties-of-O}
  Let $\vec{w} \in \sobh1(\mathcal{T}_h)^d \cap H(\mathrm{div};\Omega)$
  with $\div{\vec w}  = 0$ in $\Omega$ and $\vec w\cdot\vec n = 0$ on
  $\partial\Omega$, and let $\vec v_h \in \mathcal{P}_p(\mathcal{T}_h)^d$.
  Then the convective term satisfies
  \begin{equation}
    \mathcal{O}_h(\vec{w}; \vec v_h, \vec v_h)
    =
    \frac{\mu}{2} \sum_{e \in \mathcal{E}_h} \int_e |\vec{w} \cdot \vec{n}| 
      \ujump{\vec v_h}^2  \mathrm ds.
  \end{equation}
  Furthermore, for any $\vec{u} \in H^1_0(\Omega)^d +
  \mathcal{P}_p(\mathcal{T}_h)^d$ and any $\vec{w},\vec{w}_1$ as
  above, there exists a constant $C > 0$, independent of $h$, such
  that
  \begin{equation}
    \bigl|\mathcal{O}_h(\vec{w}; \vec{u}, \vec v_h)
      - \mathcal{O}_h(\vec{w}_1; \vec{u}, \vec v_h)\bigr|
    \leq
    C \Norm{\vec{w} - \vec{w}_1}_{1,h}
      \Norm{\vec{u}}_{1,h}
      \Norm{\vec v_h}_{1,h}.
  \end{equation}
\end{The}

\subsection{Steady-state weak formulation}

In this section we analyse existence and stability for the steady-state
version of the discrete problem. For this section only, we allow a
non-zero body force $\vec f$ in the momentum equation and define the
linear functional
\begin{equation*}
  F(\vec v_h,\vec c_h)=\langle \vec f, \vec v_h\rangle .
\end{equation*}
We introduce the nonlinear form
\begin{multline}
  \mathfrak{A}((\vec{u}_h, \vec{b}_h); (\vec{v}_h, \vec{c}_h))
  := 
  \Re \mathcal{O}_h(\vec{u}_h; \vec{u}_h, \vec{v}_h)
  +
  \beta \mathcal{A}_h(\vec{u}_h, \vec{v}_h)
  -
  \frac{1 - \beta}{2 \Wi} \left( \mathcal{O}_h(\vec{b}_h; \vec{b}_h, \vec{v}_h)
  -
  \mathcal{O}_h(\vec{b}_h; \vec{v}_h, \vec{b}_h) \right)  
  \\
  +
  \frac{1 - \beta}{2 \Wi} \left(
  \mathcal{O}_h(\vec{u}_h; \vec{b}_h, \vec{c}_h)
  -
  \frac{1}{2} \mathcal{O}_h(\vec{b}_h; \vec{u}_h, \vec{c}_h)
  +
  \frac{1}{2} \mathcal{O}_h(\vec{b}_h; \vec{c}_h, \vec{u}_h)
  \right)
  +
  \frac{1 - \beta}{\Wi^2} \ltwop{\vec b_h}{\vec c_h},
\end{multline}
and the bilinear form
\begin{equation}
  \mathfrak{B}((\vec{v}_h, \vec{c}_h), (q_h, s_h)) := 
  \mathcal{C}(\vec{v}_h, q_h) + \mathcal{D}(\vec{c}_h, s_h).
\end{equation}

The steady-state discrete problem then reads: find
$(\vec{u}_h, p_h, \vec{b}_h, r_h) \in
\mathcal{U}_h \times \mathcal{P}_h \times \mathcal{B}_h \times
\mathcal{R}_h$ such that
\begin{equation}
  \label{eq:discrete-system}
  \begin{aligned}
    \mathfrak{A}((\vec{u}_h, \vec{b}_h); (\vec{v}_h, \vec{c}_h))
    +
    \mathfrak{B}((\vec{v}_h, \vec{c}_h), (p_h, r_h)) &= F(\vec{v}_h, \vec{c}_h),
    \\
    \mathfrak{B}((\vec{u}_h, \vec{b}_h), (q_h, s_h)) &= 0,
  \end{aligned}
\end{equation}
for all $(\vec{v}_h, q_h, \vec{c}_h, s_h)\in
    \mathcal{U}_h \times \mathcal{P}_h \times \mathcal{B}_h \times \mathcal{R}_h$.

We equip the discrete spaces with the norms
\begin{equation}
  \Norm{(\vec{v}_h, \vec{c}_h)}_{X_h}^2 :=
  \Norm{\vec{v}_h}_{1,h}^2 + \Norm{\vec{c}_h}_{\leb2(\Omega)}^2,
  \qquad
  \Norm{(q_h, s_h)}_{Y_h}^2 :=
  \Norm{q_h}_{\leb2(\Omega)}^2 + \Norm{\nabla s_h}_{\leb2(\Omega)}^2.
\end{equation}

\begin{Lem}[Discrete stationary inf-sup]
  \label{lem:static-infsup}
  There exists a constant $C > 0$, independent of $h$, such that
  \begin{equation}
    \label{eq:megainf-sup}
    \sup_{(\vec{v}_h, \vec{c}_h) \in \mathcal{U}_h \times \mathcal{B}_h}
    \frac{\mathfrak{B}((\vec{v}_h, \vec{c}_h), (q_h, s_h))}{\Norm{(\vec{v}_h, \vec{c}_h)}_{X_h}}
    \geq
    C \Norm{(q_h, s_h)}_{Y_h},
    \quad \Foreach (q_h, s_h) \in \mathcal{P}_h \times \mathcal{R}_h.
  \end{equation}
  Furthermore, there exists a constant $C > 0$, independent of $h$, such that
  \begin{equation}
    \label{eq:megacoerc}
    \mathfrak{A}((\vec{v}_h, \vec{c}_h); (\vec{v}_h, \vec{c}_h))
    \geq
    C \Norm{(\vec{v}_h, \vec{c}_h)}_{X_h}^2,
    \quad \Foreach (\vec{v}_h, \vec{c}_h) \in \mathcal U_h \times \mathcal B_h.
  \end{equation}
\end{Lem}

\begin{proof}
  The inf-sup estimate \eqref{eq:megainf-sup} is a direct consequence
  of Lemma~\ref{lem:infsup} and the definitions of the norms
  $\|\cdot\|_{X_h}$ and $\|\cdot\|_{Y_h}$.

  To prove \eqref{eq:megacoerc}, let
  $(\vec{v}_h, \vec{c}_h) \in \mathcal U_h \times \mathcal B_h$. Using
  the coercivity of $\mathcal A_h$ and the positivity properties of
  $\mathcal{O}_h$ from Theorem~\ref{Theo:Properties-of-O}, together
  with the definition of $\mathfrak A$, we obtain
  \begin{equation}
    \begin{split}
      \mathfrak{A}((\vec{v}_h, \vec{c}_h); (\vec{v}_h, \vec{c}_h))
      &\geq      
      \alpha \beta \Norm{\vec v_h}_{1,h}^2
      +
      \frac{1-\beta}{\Wi^2} \Norm{\vec c_h}_{\leb{2}(\Omega)}^2
      \\
      &\qquad +
      \frac{1}{2} \sum_{e \in \mathcal{E}_h} \int_e \lvert\vec{v}_h \cdot \vec{n}\rvert
      \qp{\mu \Re  \ujump{\vec v_h}^2
        +
        (1-\beta) \Wi^{-1} \ujump{\vec c_h}^2
      }
       \mathrm d s
      \\
      &\geq
      \alpha \beta \Norm{\vec v_h}_{1,h}^2
      +
      \frac{1-\beta}{\Wi^2} \Norm{\vec c_h}_{\leb{2}(\Omega)}^2.
    \end{split}
  \end{equation}
  This shows that $\mathfrak A$ is coercive on
  $\mathcal U_h \times \mathcal B_h$ with respect to
  $\|\cdot\|_{X_h}$, and proves \eqref{eq:megacoerc}.
\end{proof}

From Lemma~\ref{lem:static-infsup}, a standard application of
Brouwer's fixed-point theorem (see, for example,
\cite[Theorem~3.1.13]{BGHRR24}) yields the steady-state existence
result stated below.

\begin{Cor}[Existence of the discrete problem]
  Under the assumptions of Lemma \ref{lem:static-infsup}, there exists
  at least one solution $(\vec{u}_h, p_h, \vec{b}_h, r_h) \in
  \mathcal{U}_h \times \mathcal{P}_h \times \mathcal{B}_h \times
  \mathcal{R}_h$ to the discrete system \eqref{eq:discrete-system}.
  In addition, there exists a constant $C>0$, depending only on the
  physical parameters of the problem and on the domain $\Omega$, such
  that
  \begin{equation}
    \Norm{(\vec{u}_h, p_h, \vec{b}_h, r_h)}_{X_h \times Y_h}
    \leq C \|F\|_{X_h^*}.
  \end{equation}
\end{Cor}
  
\begin{Rem}[Uniqueness]
  Uniqueness of the discrete solution can only be inferred if $\Re$,
  $\Wi$ and the data are sufficiently small. Such smallness
  conditions are well known already in the incompressible
  Navier-Stokes regime, where uniqueness holds at low $\Re$ but may
  fail at high $\Re$ due to the growth of nonlinear effects
  \cite{BGHRR24}. In the present coupled system, a comparable balance
  between dissipation and nonlinear growth is required to guarantee
  uniqueness.
\end{Rem}

\subsection{Spatially semidiscrete weak formulation}

Given solenoidal initial conditions
$(\vec{u}_0, \vec{b}_0) \in \mathcal{U} \times \mathcal{B}$ with
$\div \vec u_0 = 0$ and $\div \vec b_0 = 0$ in $\Omega$, we define
discrete initial data by
\[
  \vec u_{h,0} := \Pi_{\mathcal{U}_h}\vec u_0,
  \qquad
  \vec b_{h,0} := \Pi_{\mathcal{B}_h}\vec b_0,
\]
where $\Pi_{\mathcal{U}_h}$ and $\Pi_{\mathcal{B}_h}$ are the
commuting projections associated with the discrete de~Rham complex.
The spatially semidiscrete problem then reads: find
$(\vec{u}_h, p_h, \vec{b}_h, r_h)(t) \in \mathcal{U}_h \times
\mathcal{P}_h \times \mathcal{B}_h \times \mathcal{R}_h$ such that
\begin{equation}
  \label{eq:disweakform-frak}
  \begin{aligned}
    \Re \ltwop{\partial_t \vec{u}_h}{\vec{v}_h}
    + \ltwop{\partial_t \vec{b}_h}{\vec{c}_h}
    + \mathfrak{A}((\vec{u}_h, \vec{b}_h), (\vec{v}_h, \vec{c}_h))
    + \mathfrak{B}((\vec{v}_h, \vec{c}_h), (p_h, r_h))
    &= 0,
    \\
    \mathfrak{B}((\vec{u}_h, \vec{b}_h), (q_h, s_h)) &= 0,
  \end{aligned}
\end{equation}
for all test functions
$(\vec{v}_h, q_h, \vec{c}_h, s_h) \in \mathcal{U}_h \times \mathcal{P}_h
\times \mathcal{B}_h \times \mathcal{R}_h$, with initial data
$\vec u_h(0) = \vec u_{h,0}$ and $\vec b_h(0) = \vec b_{h,0}$.

\begin{The}[Existence, uniqueness and energy law for the spatially semidiscrete problem]
  \label{the:existenceSemi}
  Assume that $(\vec{u}_0, \vec{b}_0) \in \mathcal{U} \times \mathcal{B}$
  with $\div \vec u_0 = 0$ and $\div \vec b_0 = 0$ in $\Omega$, and
  let the discrete initial data be given by the projections
  $\vec u_{h,0}$ and $\vec b_{h,0}$ defined above. Then there exists a
  unique solution
  \[
    (\vec{u}_h(t), p_h(t), \vec{b}_h(t), r_h(t))
    \in
    \mathcal{U}_h \times \mathcal{P}_h \times \mathcal{B}_h \times
    \mathcal{R}_h
  \]
  to the spatially semidiscrete problem \eqref{eq:disweakform-frak} on
  $[0, T]$, depending continuously on the initial data.

  Moreover, the solution satisfies the discrete energy law
  \begin{align*}
    \ddt\qp{
      \frac{\Re}{2} \|\vec{u}_h\|^2_{L^2(\Omega)}
      +
      \frac{1-\beta}{2\Wi}\|\vec{b}_h\|^2_{L^2(\Omega)}
    }
    =&
    - \beta \mathcal A_h(\vec u_h, \vec u_h)
    -
    \frac{1-\beta}{\Wi^2} \|\vec{b}_h\|^2_{L^2(\Omega)}
    \\
    &\quad -
    \Re \mathcal{O}_h(\vec{u}_h; \vec{u}_h, \vec{u}_h)
    -
    \frac{1-\beta}{\Wi} \mathcal{O}_h\qp{\vec{u}_h; \vec{b}_h, \vec{b}_h},
  \end{align*}
  so that, in particular, the discrete total energy is non-increasing
  in time.
\end{The}
  
\begin{Rem}[Discrete dissipation]
  Comparing the discrete energy law above with the continuous energy
  law \eqref{eq:energy}, we observe the additional dissipative
  contributions introduced by the upwind discretisation of the
  convective terms. In particular, the facet terms in
  $\mathcal{O}_h(\vec u_h;\cdot,\cdot)$ provide extra damping through
  jump penalisation, which is absent at the continuous level.
\end{Rem}

\begin{Proof}[of Theorem \ref{the:existenceSemi}]
  To derive the discrete energy law, in \eqref{eq:disweakform-frak} we
  take $\vec{v}_h = \vec{u}_h$ in the first equation and $\vec{c}_h =
  \tfrac{1-\beta}{\Wi} \vec{b}_h$ in the second one, and add the
  resulting identities. Choosing $q_h = p_h$ and $s_h = r_h$, and
  using the cancellation of the skew-symmetric convective terms
  together with the coercivity of $\mathcal A_h$ and the properties of
  $\mathcal O_h$ from Theorem~\ref{Theo:Properties-of-O}, yields the
  stated discrete energy law.

  For well-posedness, we follow the approach in, for example,
  \cite[Section~5]{BF07}. We first introduce the discrete kernel
  \begin{equation}
    \mathcal{Z}_h
    :=
    \bigl\{
      (\vec v_h,\vec c_h)\in\mathcal{U}_h\times\mathcal{B}_h :
      \mathfrak{B}((\vec{v}_h, \vec{c}_h), (q_h, s_h))=0
      \ \forall (q_h, s_h)\in \mathcal{P}_h\times\mathcal{R}_h
    \bigr\}.
  \end{equation}
  The semidiscrete weak formulation \eqref{eq:disweakform-frak} is
  equivalent to an abstract system of ODEs posed on $\mathcal{Z}_h$:
  find $\vec{z}_h(t)\in\mathcal{Z}_h$ such that
  $\vec z_h(0)= (\vec u_{h,0}, \vec b_{h,0})$ and
  \begin{equation}
    \vec M \frac{\mathrm d\vec{z}_h}{\mathrm d t} + \vec A(\vec{z}_h) = 0,
  \end{equation}
  where $\vec M$ is the (symmetric) positive definite mass matrix
  associated with the $L^2$ inner product, and $\vec A(\vec{z}_h)$ is
  the nonlinear operator induced by $\mathfrak{A}$ restricted to
  $\mathcal{Z}_h$.

  Using Theorem~\ref{Theo:Properties-of-O} and the boundedness of the
  remaining terms in $\mathfrak A$, we infer that $\vec A$ is
  Lipschitz continuous on $\mathcal{Z}_h$ with respect to the
  $\|\cdot\|_{X_h}$-norm, that is,
  \begin{equation}
    \Norm{\vec A(\vec{z}_h) - \vec A(\vec{z}_h')}_{X_h^*}
    \leq
    L \Norm{\vec{z}_h - \vec{z}_h'}_{X_h},
  \end{equation}
  for some constant $L>0$ independent of $h$. The Picard-Lindel\"of
  theorem then guarantees the existence and uniqueness of a solution
  $\vec{z}_h(t)$ on $[0,T]$, and therefore of
  $(\vec{u}_h(t), p_h(t), \vec{b}_h(t), r_h(t))$, completing the
  proof.
\end{Proof}

\section{An IMEX temporal discretisation}
\label{sec:fulldis}

The choice of temporal discretisation impacts the properties of the
fully discrete energy law. Structure-preserving schemes, such as the
Crank-Nicolson method, can be used to retain a discrete analogue of
the continuum energy balance. On the other hand, dissipative schemes
such as backward differentiation formulae introduce additional
numerical diffusion. This added dissipation can be beneficial in flows
where excessive energy growth from nonlinearities, such as elastic
turbulence, must be controlled. However, these methods are inherently
nonlinear, requiring iterative solvers at each time step and thus
increasing computational cost.

In this section we introduce an implicit-explicit (IMEX) method that
preserves a discrete energy law while keeping the scheme linear at
each time step, offering a balance between stability and computational
efficiency.

To discretise the system in time, let $N \in \mathbb{N}$ be the number
of time steps, and define temporal nodes
$0 = t^0 < t^1 < \dots < t^N = T$. For simplicity we take a uniform
time step $\tau := t^{n} - t^{n-1}$, and for a generic function $u$
write $u^n := u(t^n)$ for its value at time $t^n$.

\subsection{Time-stepping scheme}

Given discrete initial conditions $\vec u_h^0 \in \mathcal{U}_h$ and
$\vec b_h^0 \in \mathcal{B}_h$, for each time step
$n = 0, 1, 2, \dots, N-1$ we seek
$(\vec u_h^{n+1}, p_h^{n+1}, \vec b_h^{n+1}, r_h^{n+1}) \in
\mathcal{U}_h \times \mathcal{P}_h \times \mathcal{B}_h \times
\mathcal{R}_h$ such that
\begin{equation}
  \label{eq:fulldis}
  \begin{split}
    \Re \ltwop{\frac{\vec u_h^{n+1} - \vec u_h^n}{\tau}}{\vec v_h}
    &+ \Re \mathcal{O}_h(\vec u_h^n; \vec u_h^{n+1}, \vec v_h)
    + \mathcal{C}(\vec v_h, p_h^{n+1}) \\
    &=
    - \beta \mathcal{A}_h(\vec u_h^{n+1}, \vec v_h)
    + \frac{1 - \beta}{2\Wi} \mathcal{O}_h(\vec b_h^n; \vec b_h^{n+1}, \vec v_h)
    - \frac{1 - \beta}{2\Wi} \mathcal{O}_h(\vec b_h^n; \vec v_h, \vec b_h^{n+1}),\\
    \ltwop{\frac{\vec b_h^{n+1} - \vec b_h^n}{\tau}}{\vec c_h}
    &+ \mathcal{O}_h(\vec u_h^n; \vec b_h^{n+1}, \vec c_h)
    - \frac{1}{2} \mathcal{O}_h(\vec b_h^n; \vec u_h^{n+1}, \vec c_h)
    + \frac{1}{2} \mathcal{O}_h(\vec b_h^n; \vec c_h, \vec u_h^{n+1}) \\
    &\quad + \mathcal{D}(\vec c_h, r_h^{n+1})
    + \frac{1}{\Wi} \ltwop{\vec b_h^{n+1}}{\vec c_h}
    = 0, \\
    \mathcal{C}(\vec u_h^{n+1}, q_h) &= 0, \\
    \mathcal{D}(\vec b_h^{n+1}, s_h) &= 0,
  \end{split}
\end{equation}
for all test functions $(\vec v_h, q_h, \vec c_h, s_h) \in
\mathcal{U}_h \times \mathcal{P}_h \times \mathcal{B}_h \times
\mathcal{R}_h$.

\begin{Rem}[Consistency]
  The IMEX scheme \eqref{eq:fulldis} is first-order accurate in time
  due to the explicit treatment of the nonlinear terms. Higher-order
  temporal accuracy may be achieved by replacing the Euler
  time-stepping with a higher-order IMEX Runge-Kutta method, as
  described in \cite{ascher1995implicit}.
\end{Rem}

\begin{The}[Existence, uniqueness and discrete energy law]
  \label{thm:discrete-well-posedness-energy}
  Let $(\vec{u}^0, \vec{b}^0) \in \mathcal{U} \times \mathcal{B}$ be
  initial conditions satisfying the divergence-free constraints, and
  let the corresponding discrete initial data be
  $\vec u_h^0 = \Pi_{\mathcal U_h}\vec u^0$ and
  $\vec b_h^0 = \Pi_{\mathcal B_h}\vec b^0$. Then, for each
  $n = 0, 1, \dots, N-1$ there exists a unique solution
  $(\vec u_h^{n+1}, p_h^{n+1}, \vec b_h^{n+1}, r_h^{n+1}) \in
  \mathcal{U}_h \times \mathcal{P}_h \times \mathcal{B}_h \times
  \mathcal{R}_h$ to the fully discrete system \eqref{eq:fulldis}.

  Furthermore, the solution satisfies the discrete energy law
  \begin{multline}
  \label{eq:discrete_energy_law}
    \frac{\Re}{2} \Norm{\vec u_h^{n+1}}^2_{\leb 2(\Omega)}
    +
    \frac{1-\beta}{2\Wi} \Norm{\vec b_h^{n+1}}^2_{\leb 2(\Omega)} 
    =
    \frac{\Re}{2} \Norm{\vec u_h^0}^2_{\leb 2(\Omega)}
    +
    \frac{1-\beta}{2\Wi} \Norm{\vec b_h^0}^2_{\leb 2(\Omega)} \\
    - \sum_{j=0}^n 
    \Bigg(
      \frac{\Re}{2} \Norm{\vec{u}_h^{j+1} - \vec{u}_h^j}_{\leb 2(\Omega)}^2
      +
      \frac{1-\beta}{2\Wi} \Norm{\vec{b}_h^{j+1} - \vec{b}_h^j}_{\leb 2(\Omega)}^2
    \Bigg) \\
    - \tau 
    \sum_{j=0}^n 
    \Bigg(
      \beta \mathcal{A}_h(\vec u_h^{j+1}, \vec u_h^{j+1})
      +
      \frac{1-\beta}{\Wi^2} \Norm{\vec b_h^{j+1}}^2_{\leb 2(\Omega)} \\
      +
      \sum_{e \in \mathcal{E}_h} 
      \int_e 
      \Big(
        \mu \Re \lvert\vec{u}_h^j \cdot \vec{n}\rvert  \ujump{\vec{u}_h^{j+1}}^2
        +
        \frac{1-\beta}{\Wi} \lvert\vec{u}_h^j \cdot \vec{n}\rvert  \ujump{\vec{b}_h^{j+1}}^2
      \Big)  \mathrm d s
    \Bigg).
  \end{multline}
\end{The}

\begin{Proof}
  We first derive the discrete energy law. In \eqref{eq:fulldis}, take
  $\vec{v}_h = \vec{u}_h^{n+1}$ in the first equation and
  $\vec{c}_h = \tfrac{1-\beta}{\Wi} \vec{b}_h^{n+1}$ in the second
  equation, and add the resulting identities. Using $\mathcal C$ and
  $\mathcal D$ with $(q_h,s_h) = (p_h^{n+1},r_h^{n+1})$ gives
  \begin{multline}
    \Re \ltwop{\frac{\vec{u}_h^{n+1} - \vec{u}_h^n}{\tau}}{\vec{u}_h^{n+1}}
    + \frac{1-\beta}{\Wi}\ltwop{\frac{\vec{b}_h^{n+1} - \vec{b}_h^n}{\tau}}{\vec{b}_h^{n+1}}
    + \beta \mathcal{A}_h(\vec{u}_h^{n+1}, \vec{u}_h^{n+1})
    \\
    + \frac{1-\beta}{\Wi^2} \Norm{\vec{b}_h^{n+1}}^2_{\leb 2(\Omega)} 
    + \sum_{e \in \mathcal{E}_h} \int_e 
      \Big(
        \mu \Re \lvert\vec{u}_h^n \cdot \vec{n}\rvert  \ujump{\vec{u}_h^{n+1}}^2 
        +
        \frac{1-\beta}{\Wi}\lvert\vec{u}_h^n \cdot \vec{n}\rvert  \ujump{\vec{b}_h^{n+1}}^2
      \Big) \mathrm d s
    = 0.  
  \end{multline}
  The time-derivative terms are rewritten using
  \[
    (x-y) x = \tfrac12 x^2 - \tfrac12 y^2 + \tfrac12 (x-y)^2,
  \]
  which yields
  \begin{multline}
    \frac{\Re}{2\tau} \Bigl( \Norm{\vec{u}_h^{n+1}}^2_{\leb 2(\Omega)}
    -
    \Norm{\vec{u}_h^n}^2_{\leb 2(\Omega)} \Bigr)
    +
    \frac{1-\beta}{2\tau\Wi} \Bigl(
      \Norm{\vec{b}_h^{n+1}}^2_{\leb 2(\Omega)}
      -
      \Norm{\vec{b}_h^n}^2_{\leb 2(\Omega)}
    \Bigr)
    \\
    +
    \frac{\Re}{2}\Norm{\vec{u}_h^{n+1} - \vec{u}_h^n}_{\leb 2(\Omega)}^2
    +
    \frac{1-\beta}{2\Wi}\Norm{\vec{b}_h^{n+1} - \vec{b}_h^n}_{\leb 2(\Omega)}^2
    \\
    +
    \beta \mathcal{A}_h(\vec{u}_h^{n+1}, \vec{u}_h^{n+1})
    + \frac{1-\beta}{\Wi^2} \Norm{\vec{b}_h^{n+1}}^2_{\leb 2(\Omega)} 
    + \sum_{e \in \mathcal{E}_h} \int_e 
      \Big(
        \mu \Re \lvert\vec{u}_h^n \cdot \vec{n}\rvert  \ujump{\vec{u}_h^{n+1}}^2 
        +
        \frac{1-\beta}{\Wi} \lvert\vec{u}_h^n \cdot \vec{n}\rvert  \ujump{\vec{b}_h^{n+1}}^2
      \Big) \mathrm d s
    = 0.
  \end{multline}
  Multiplying by $\tau$ and rearranging gives a one-step energy
  balance. Summing this identity from $j=0$ to $j=n$ and telescoping
  the leading terms yields \eqref{eq:discrete_energy_law}.

  To establish existence and uniqueness at each time step, we observe
  that \eqref{eq:fulldis} defines, for fixed data
  $(\vec u_h^n,\vec b_h^n)$, a linear system for
  $(\vec u_h^{n+1}, p_h^{n+1}, \vec b_h^{n+1}, r_h^{n+1})$. The
  discrete energy identity implies that if the right-hand side
  corresponding to the previous time step is homogeneous, then the
  only solution is the trivial one, so the associated linear operator
  is injective. Since the operator acts on a finite-dimensional space,
  injectivity implies bijectivity, and thus a unique solution exists
  for each time step. By induction on $n$, this gives existence and
  uniqueness of the fully discrete solution sequence.
\end{Proof}

\section{Numerical Experiments}
\label{sec:numerics}

In this section we numerically investigate the IMEX discretisation
\eqref{eq:fulldis}, implemented in the FEniCS finite element
software. In Section \ref{oldroyd-sec:benchmark} we verify
convergence, and in Sections
\ref{oldroyd-sec:lid-driven-cavity}-\ref{oldroyd-sec:contraction-flow}
we demonstrate the stability of the method in challenging
physically-inspired experiments for a range of Weissenberg numbers. We
focus on low Reynolds number regimes and examine the impact of varying
the elasticity ratio $\beta$ and the Weissenberg number $\Wi$.

For each of the physically-inspired experiments we first consider the
discretisation \eqref{eq:fulldis} applied to the reduced system
\eqref{eq:reduced_model_divfree} as derived in
Section~\ref{sec:modelreduction}. We then repeat the experiment with
an additional explicit term
\[
  \frac{1-\gamma}{\Wi}\,\frac{\vec{b}_h^n}{\abs{\vec{b}_h^n}}
\]
on the right-hand side of the equation for $\vec{b}_h^{n+1}$ (with the
convention that this contribution is set to zero when
$\vec b_h^n = \vec 0$), where $\gamma = 10^{-6}$. In the full
Oldroyd--B model the relaxation term $-\Wi^{-1}(\sig - \vec I)$ drives
the conformation towards the isotropic state $\sig = \vec I$, so that
at rest the principal eigenvalue tends to $1$. In the strong-stretch,
uniaxial reduction used here, where
$\sig \approx \vec b\otimes\vec b$ and the identity contribution has
been neglected, the reduced equation contains only the damping
$-\Wi^{-1}\vec b$ and therefore has an equilibrium at $\vec b = \vec
0$, corresponding to vanishing polymeric stress and an effectively
Newtonian response. The additional explicit term can be viewed as a
simple surrogate for the omitted identity contribution: in the absence
of flow it tends to stabilise the amplitude $\abs{\vec b}$ near a
non-zero value of order one (roughly $\abs{\vec b}\approx 1-\gamma$),
while leaving the directional dynamics governed by transport and
stretching. Our numerical results show that including this term yields
flow patterns that are closer to those observed in full
Oldroyd--B simulations. The term is treated explicitly so as not to
change the discrete energy law established for the unmodified scheme,
and the factor $1-\gamma<1$ introduces a small net damping that
prevents unbounded growth of the polymeric stretch.

We take $\mu = 1$ throughout and, unless otherwise specified, all
examples are conducted in two spatial dimensions using unstructured
Delaunay triangular meshes. The linear system at each time step is
solved with a direct solver.

\subsection{Convergence to a smooth solution on a uniform mesh}
\label{oldroyd-sec:benchmark}

\newcommand{\upstreamlength}{L_1}
\newcommand{\contractionlength}{L_2}
\newcommand{\upstreamwidth}{H_1}
\newcommand{\contractionwidth}{H_2}
\newcommand{\contractionratio}{CR}
\newcommand{\inflowvelocity}{u_-}
\newcommand{\outflowvelocity}{u_+}

We first verify the accuracy of the spatial discretisation by studying
the error with respect to the manufactured solution
\begin{align}
  \vec{u}(\vec x, t) &= \Transpose{\qp{1,1}} +\sin{\frac{\pi t}{2}}\Transpose{\left(-\cos{\pi x}\sin{\pi y}, \sin{\pi x}\cos{\pi y}\right)},
  &\vec{b}(\vec x, t) &= \sin{\frac{\pi t}{2}}\Transpose{\left(-x, y\right)}, \\
  p(\vec x,t) &= \sin{\frac{\pi t}{2}}\qp{x + y - 1},
  &r(\vec x,t) &= 0,
\end{align}
which yields divergence-free fields $\vec{u}$ and $\vec{b}$. The
spatial domain is the unit square and the final time is $T = 1$. We
take $\Re = \Wi = 1$ and $\beta = 1/2$. The flow is initialised with
$\vec{u}_0 = \vec 0$ and $\vec{b}_0 = \Transpose{\qp{1, 1}}$. We take
the interior penalty parameter $\sigma_{\mathrm{i}} = 10$ and
$\sigma_{\mathrm{b}} = 10^5$.

\begin{figure}[h]
  \centering
  \begin{tikzpicture}
    \node [draw,fill=white] at (0.8,0.3) {\shortstack[l]{
        \ref{p1} $k = 1$ \quad
        \ref{p2} $k = 2$}};
  \end{tikzpicture}
  \vspace{1em}
  
  \subcaptionbox{}
  {\begin{tikzpicture}[scale=0.9]
  \begin{axis}[
      cycle list/Dark2,
      thick,
      xmode=log,
      ymode=log,
      xlabel=$h$,
      ylabel=$\Norm{\vec u - \vec u_h}_{\mathcal{U}}$,
      grid=both,
      minor grid style={gray!25},
      major grid style={gray!25},
      legend style=none,
    ]
    \addplot[color=pone, mark=triangle] table[x expr={1.0/\thisrow{h}}, y=u energy, col sep=comma] {data/errors_p1.csv};
    \addlegendimage{only marks, color=pone, mark=triangle, thick}
    \label{p1}
    \addplot[color=ptwo, mark=square] table[x expr={1.0/\thisrow{h}}, y=u energy, col sep=comma] {data/errors_p2.csv};
    \addlegendimage{only marks, color=ptwo, mark=square, thick}
    \label{p2}

    \ConvergenceTriangle{10}{5e-1}{2.5}{1};
    \ConvergenceTriangle{13}{2.5e-2}{2}{2};
  \end{axis}
\end{tikzpicture}

}
  \subcaptionbox{\label{oldroyd-fig:b-convergence}}
  {\begin{tikzpicture}[scale=0.9]
  \begin{axis}[
      cycle list/Dark2,
      thick,
      xmode=log,
      ymode=log,
      xlabel=$h$,
      ylabel=$\Norm{\vec b - \vec b_h}_{\mathcal{B}}$,
      grid=both,
      minor grid style={gray!25},
      major grid style={gray!25},
      legend style=none,
    ]
    \addplot[color=pone, mark=triangle] table[x expr={1.0/\thisrow{h}}, y=B energy, col sep=comma] {data/errors_p1.csv};
    \addlegendimage{only marks, color=pone, mark=triangle, thick}
    \label{p1}
    \addplot[color=ptwo, mark=square] table[x expr={1.0/\thisrow{h}}, y=B energy, col sep=comma] {data/errors_p2.csv};
    \addlegendimage{only marks, color=ptwo, mark=square, thick}
    \label{p2}

    \ConvergenceTriangle{10}{3e-2}{2.5}{1.5};
    \ConvergenceTriangle{12}{3e-3}{2}{2};
  \end{axis}
\end{tikzpicture}

}
  \caption{Calculated errors for (a) $\vec u_h$ and (b) $\vec b_h$,
  using finite elements of degree $k = 1,2$, for the manufactured
  solution test of Section \ref{oldroyd-sec:benchmark}.}
  \label{oldroyd-fig:benchmark}
\end{figure}
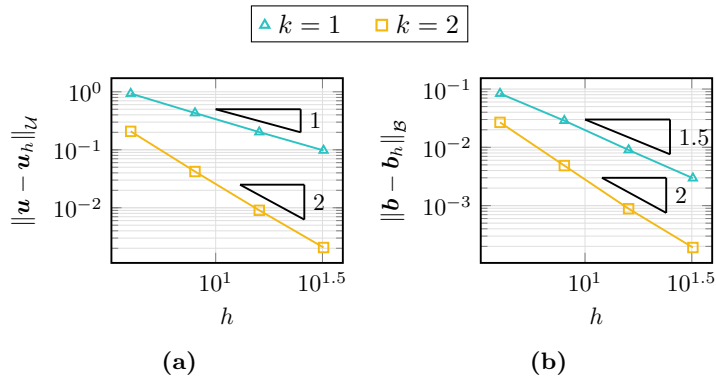

For elements of polynomial degree $k = 1,2$, the discrete problem
\eqref{eq:fulldis} is solved, with appropriate additional forcing
terms, on a sequence of uniformly refined meshes with mesh parameters
$h^{-1} = 4, 8, 16, 32$. The time step is coupled to the mesh as
$\tau = h^{k+1}$. The resulting errors in $\vec{u}_h$ and $\vec{b}_h$
are measured in the norms
\begin{align}
  \Norm{\vec{w}}_{\mathcal{U}}^2
  &:= \Norm{\vec{w}}_{L^\infty(L^2(\Omega))}^2
    + \frac{\beta}{\Re}\Norm{\nabla\vec{w}}_{L^2(L^2(\Omega))}^2, \\
  \Norm{\vec{w}}_{\mathcal{B}}^2
  &:= \Norm{\vec{w}}_{L^\infty(L^2(\Omega))}^2
    + \frac{1-\beta}{\Wi^2}\Norm{\vec{w}}_{L^2(L^2(\Omega))}^2,
\end{align}
and are plotted in Figure \ref{oldroyd-fig:benchmark}.

The results demonstrate convergence at a rate
$\mathcal{O} \left(h^k\right)$ for
$\Norm{\vec u - \vec u_h}_{\mathcal{U}}^2$, and at a rate between
$\mathcal{O} \left(h^k\right)$ and
$\mathcal{O} \left(h^{k+1}\right)$ for
$\Norm{\vec b - \vec b_h}_{\mathcal{B}}^2$. This rate is expected to
depend on the parameter choices. For instance, we have chosen
$\beta = 1/2$ to give a balance of inertial and elastic dynamics.
Moving towards a purely elastic regime as $\beta \to 0$, or towards
Newtonian flow as $\beta \to 1$, changes the character of the problem
and hence its convergence properties.

\subsection{Flow inside a lid-driven cavity}
\label{oldroyd-sec:lid-driven-cavity}

In the first physically motivated experiment we consider the
lid-driven cavity problem, in which an initially stationary fluid is
forced into motion by the movement of the upper boundary of a
rectangular domain. In the context of viscoelastic fluids this problem
has been widely studied, see for example
\cite{pan2009simulation,ashby2025discretisation} and the references
therein. It is a challenging test case due to the emergence of a
boundary layer along the moving boundary, particularly near the
corners, where large polymeric stresses accumulate and can trigger
instabilities.

\begin{figure}[h]
  \begin{tikzpicture}
  \def\length{3.0};
  \def\height{3.0};

  \colorlet{electrodeonecolour}{blue}
  \colorlet{electrodetwocolour}{red}
  \colorlet{lidcolour}{orange}
  \colorlet{noslipcolour}{black}
  
  \coordinate (electrode-one-corner-bottom) at (0.0, 0.0);
  \draw[noslipcolour] (electrode-one-corner-bottom) -- ([shift={(0.0, \height)}]electrode-one-corner-bottom) coordinate (electrode-one-corner-top);
  \draw[lidcolour, very thick] (electrode-one-corner-top) -- ([shift={(\length, 0.0)}]electrode-one-corner-top) coordinate (electrode-two-corner-top) node[black, midway, above] {$\Gamma_D$};
  \draw[noslipcolour] (electrode-two-corner-top) -- ([shift={(0.0, -\height)}]electrode-two-corner-top) coordinate (electrode-two-corner-bottom);
  \draw[noslipcolour] (electrode-two-corner-bottom) -- (electrode-one-corner-bottom) node[black, midway, below] {$\Gamma_0$};

    \fill[lidcolour, opacity=0.2] (electrode-one-corner-top) rectangle ([shift={(0, 0.1)}]electrode-two-corner-top);
    
  \draw plot[smooth cycle, -latex, tension=0.9] coordinates {(0.5 * \length, 0.9 * \height) (0.85 * \length, 0.7 * \height) (0.5 * \length, 0.25 * \height) (0.15 * \length, 0.7 * \height) } [arrow inside={end=latex, opt={scale=1.5}} {0.2, 0.7}];
  \draw plot[smooth cycle, -latex, tension=1.0] coordinates {(0.5 * \length, 0.82 * \height) (0.75 * \length, 0.68 * \height) (0.5 * \length, 0.4 * \height) (0.25 * \length, 0.68 * \height) } [arrow inside={end=latex, opt={scale=1.5}} {0.4, 0.9}];
  \draw plot[smooth cycle, -latex, tension=1.0] coordinates {(0.5 * \length, 0.75 * \height) (0.65 * \length, 0.65 * \height) (0.5 * \length, 0.55 * \height) (0.35 * \length, 0.65 * \height) } [arrow inside={end=latex, opt={scale=1.5}} {0.3, 0.8}];
  
\end{tikzpicture}

  \caption{The lid-driven cavity problem, in which an initially
    stationary fluid is forced into motion by the movement of the
    upper boundary $\Gamma_{\vec{u}}$ of a rectangular domain. The
    remaining boundary $\Gamma_0$ is fixed.}
  \label{fig:lid-driven-cavity-setup}
\end{figure}
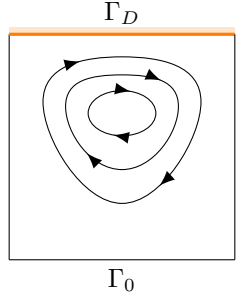

A diagram of the problem setup is shown in
Figure~\ref{fig:lid-driven-cavity-setup}. We solve \eqref{eq:fulldis}
on the unit square with polynomial degree $k = 1$, time step $\tau =
3\times10^{-3}$ and final time $T = 30$. To resolve the boundary layer
at the top of the cavity we employ a graded mesh ($341{,}622$ total
degrees of freedom), with increased resolution towards the upper edge
and especially near the upper corners. The smallest and largest
element diameters are $h_K \approx 3\times10^{-3}$ and $h_K \approx
10^{-2}$, respectively. We take the interior penalty parameter
$\sigma_{\mathrm{i}} = 10^2$ and $\sigma_{\mathrm{b}} = 10^3$. We
simulate creeping flow by choosing $\Re = 10^{-2}$, and consider two
Weissenberg regimes, $\Wi = 1/2$ and $\Wi = 1$. Unless otherwise
stated we take $\beta = 1/2$.

\begin{figure}[h!]
  \centering
  \subcaptionbox{$t^n = 2.1$}
  {\begin{subfigure}{0.2\textwidth}
      \imagewithcentreline[width=\textwidth]{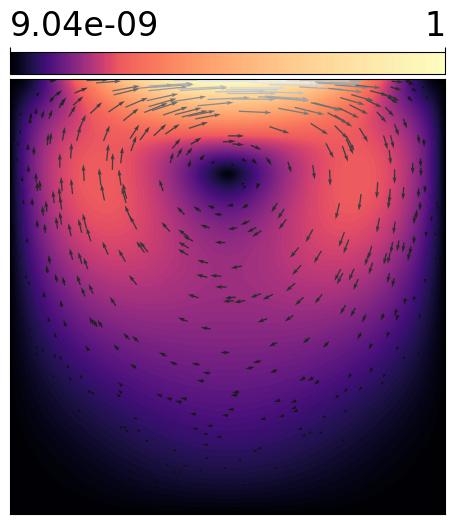} \\
      \imagewithcentreline[width=\textwidth]{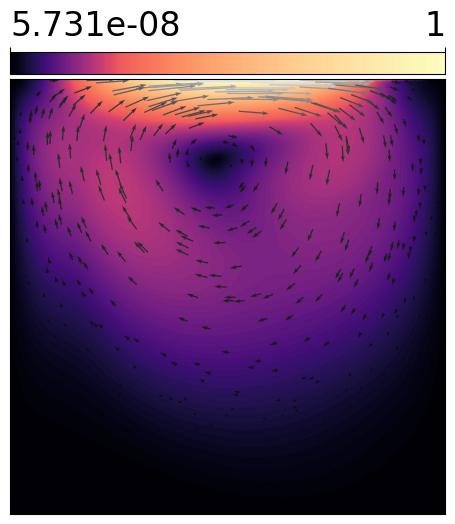}
      \imagewithcentreline[width=\textwidth]{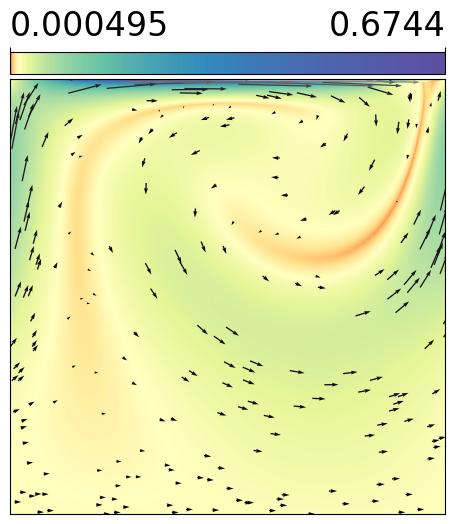} \\
      \imagewithcentreline[width=\textwidth]{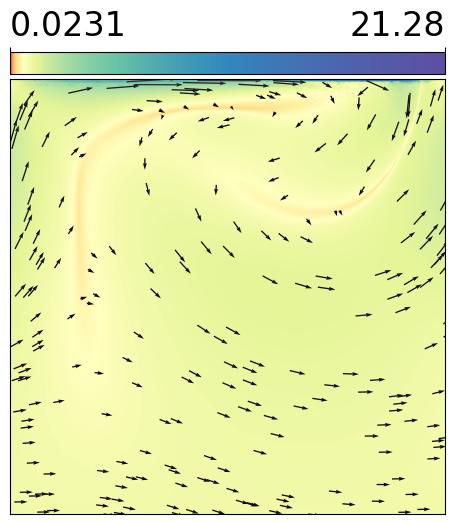}
    \end{subfigure}}
  \subcaptionbox{$t^n = 5.25$}
  {\begin{subfigure}{0.2\textwidth}
      \imagewithcentreline[width=\textwidth]{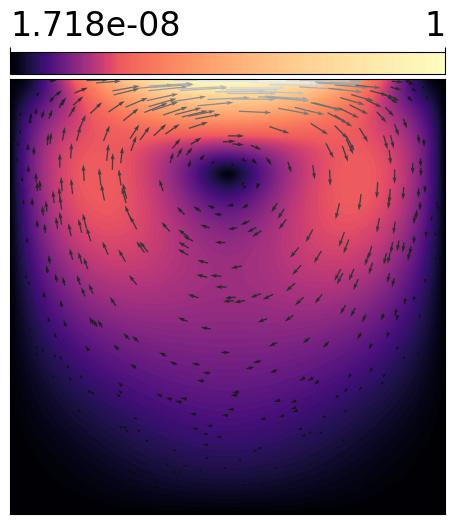} \\
      \imagewithcentreline[width=\textwidth]{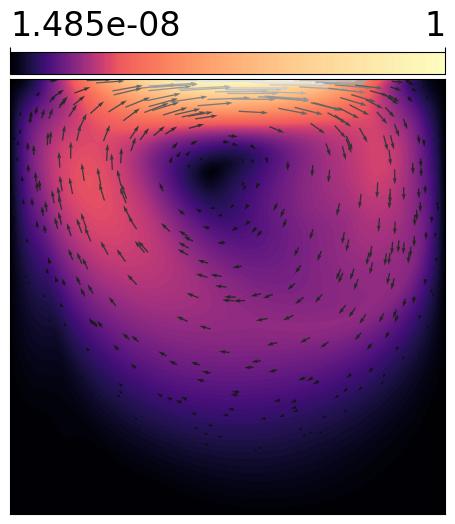}
      \imagewithcentreline[width=\textwidth]{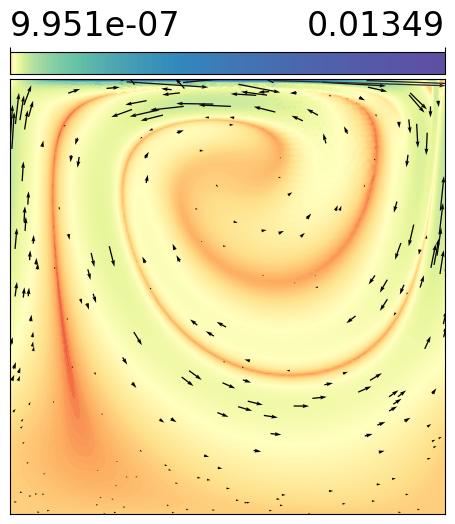} \\
      \imagewithcentreline[width=\textwidth]{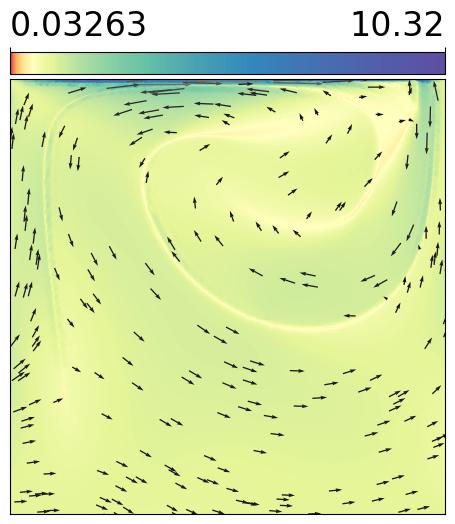}
    \end{subfigure}}
  \subcaptionbox{$t^n = 15$}
  {\begin{subfigure}{0.2\textwidth}
      \imagewithcentreline[width=\textwidth]{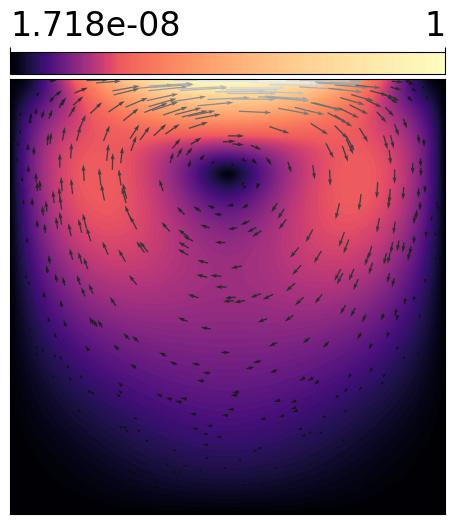} \\
      \imagewithcentreline[width=\textwidth]{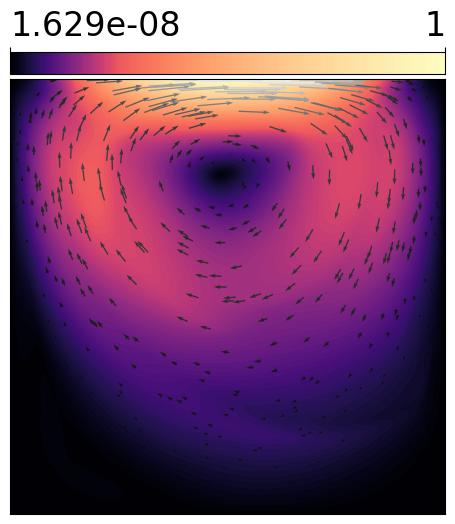}
      \imagewithcentreline[width=\textwidth]{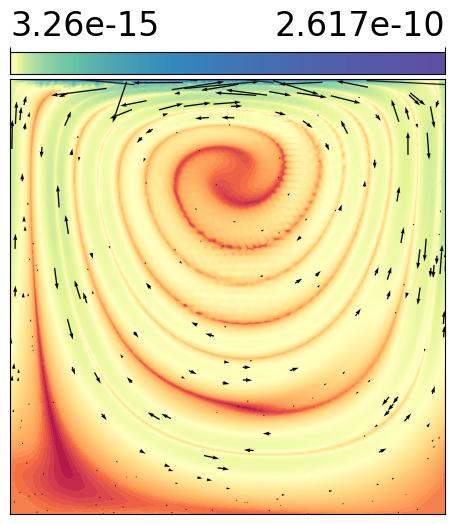} \\
      \imagewithcentreline[width=\textwidth]{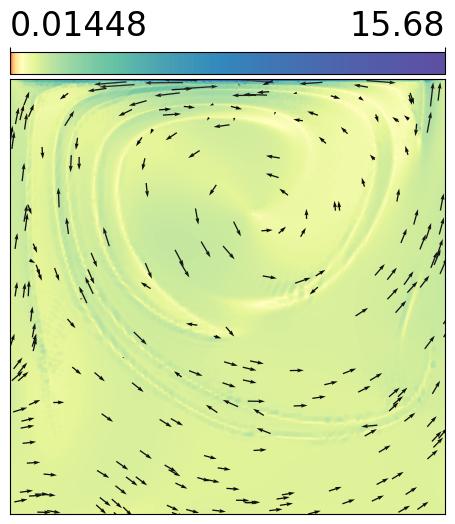}
    \end{subfigure}}
  \subcaptionbox{$t^n = 24.9$}
  {\begin{subfigure}{0.2\textwidth}
      \imagewithcentreline[width=\textwidth]{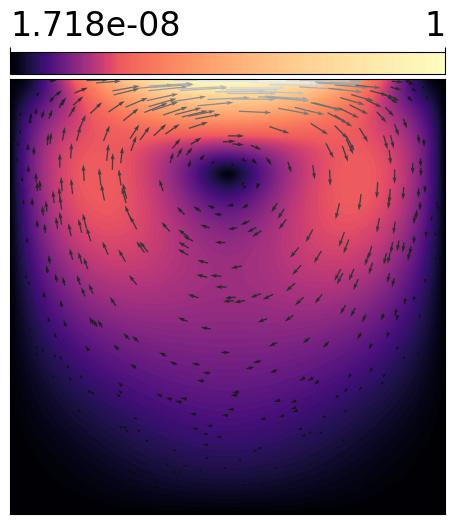} \\
      \imagewithcentreline[width=\textwidth]{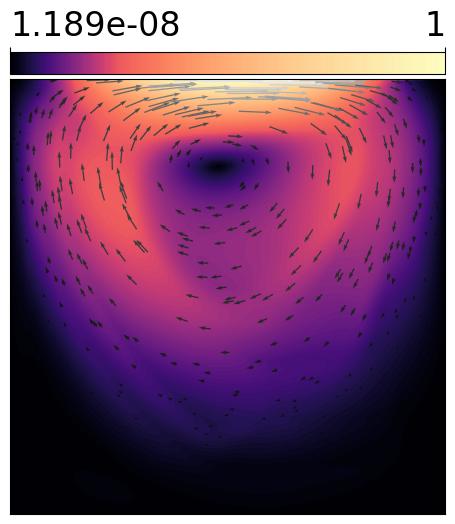}
      \imagewithcentreline[width=\textwidth]{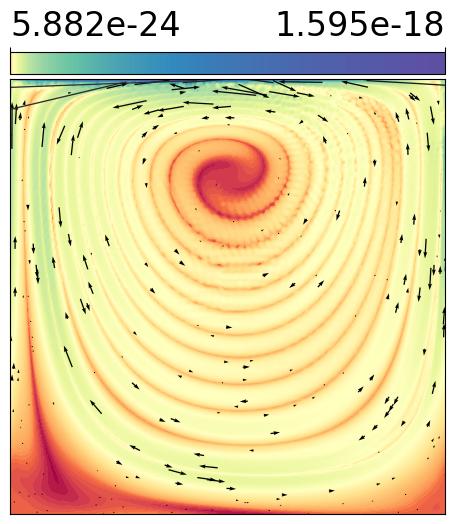} \\
      \imagewithcentreline[width=\textwidth]{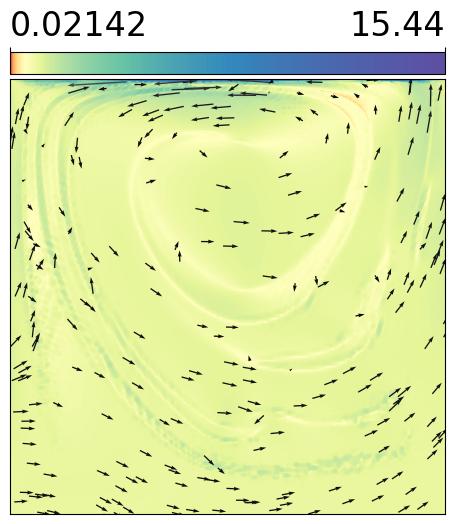}
    \end{subfigure}}
  \caption{Snapshots of $\vec{u}_h$ and $\vec{b}_h$ from the lid-driven cavity experiment in Section \ref{oldroyd-sec:lid-driven-cavity}, with $\Wi = 1/2$. The top and second rows show $\vec{u}_h$ without and with the extra forcing term, respectively, and the third and bottom rows show $\vec{b}_h$ without and with the extra forcing term, respectively. We draw particular attention to the logarithmic colour scales on each plot, which are used due to the nature of the boundary layer that arises along the top of the domain.}
  \label{oldroyd-fig:lid-driven-cavity-snapshots-small-Wi}
\end{figure}

\begin{figure}[h!]
  \centering
  \subcaptionbox{$t^n = 2.1$}
  {\begin{subfigure}{0.2\textwidth}
      \imagewithcentreline[width=\textwidth]{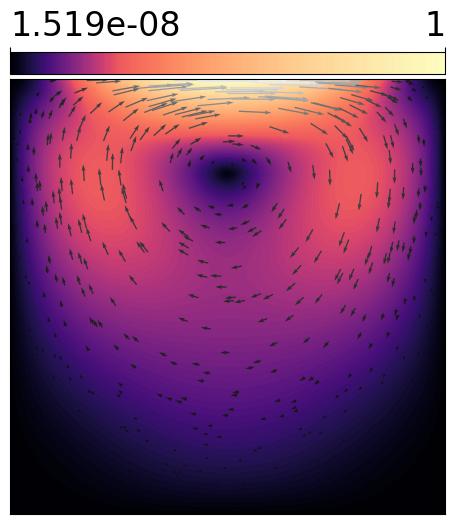} \\
      \imagewithcentreline[width=\textwidth]{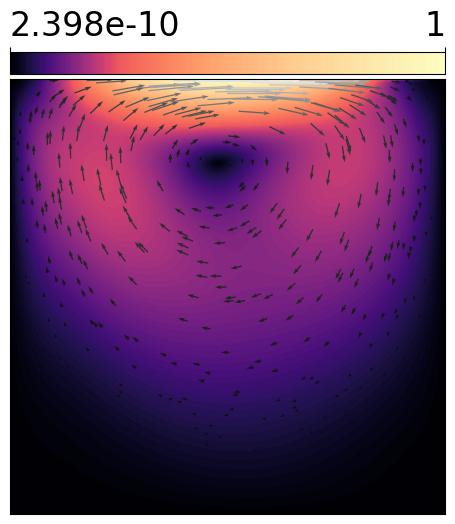}
      \imagewithcentreline[width=\textwidth]{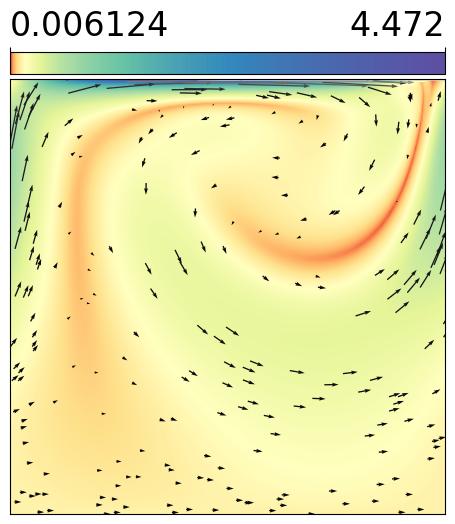} \\
      \imagewithcentreline[width=\textwidth]{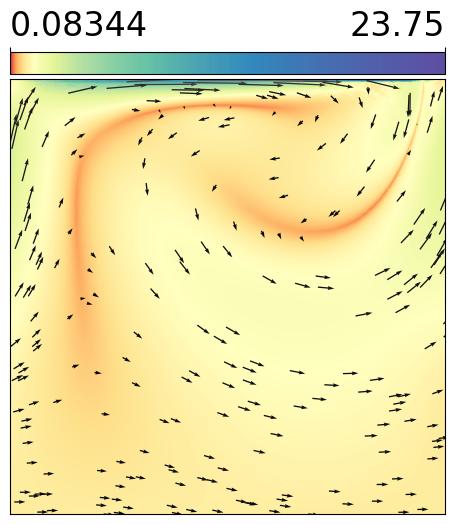}
    \end{subfigure}}
  \subcaptionbox{$t^n = 5.25$}
  {\begin{subfigure}{0.2\textwidth}
      \imagewithcentreline[width=\textwidth]{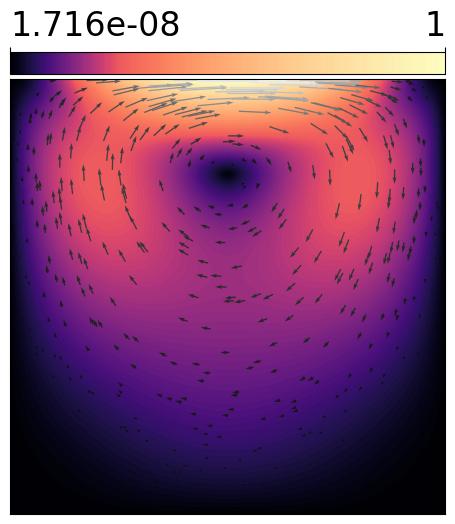} \\
      \imagewithcentreline[width=\textwidth]{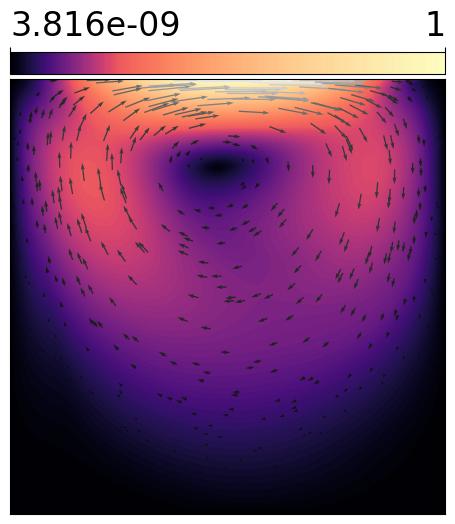}
      \imagewithcentreline[width=\textwidth]{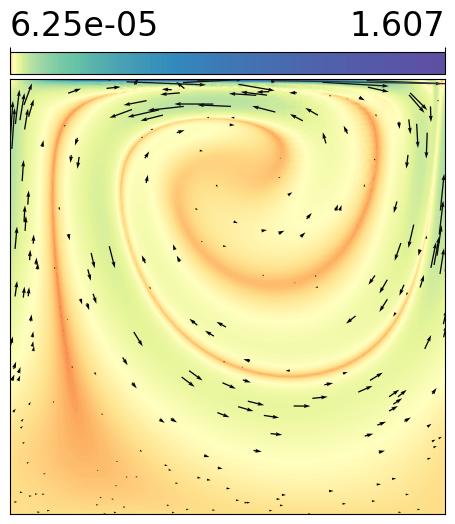} \\
      \imagewithcentreline[width=\textwidth]{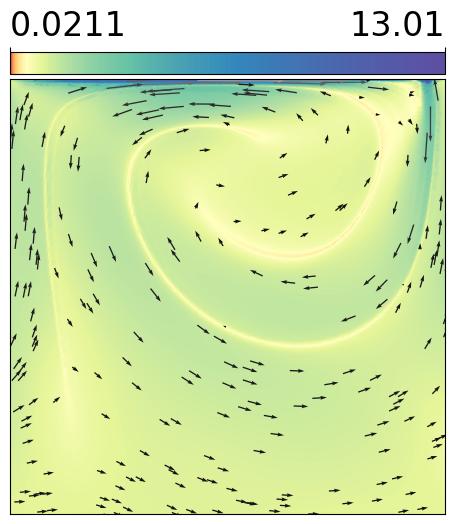}
    \end{subfigure}}
  \subcaptionbox{$t^n = 15$}
  {\begin{subfigure}{0.2\textwidth}
      \imagewithcentreline[width=\textwidth]{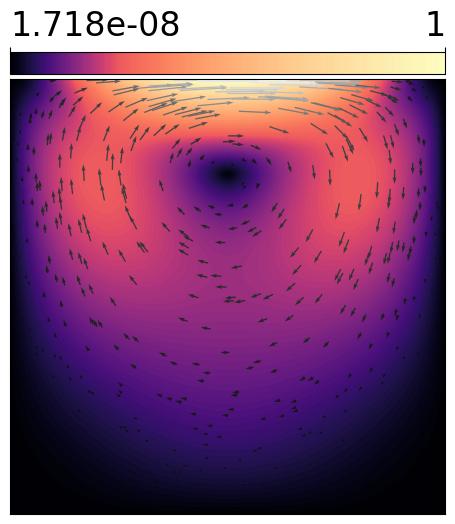} \\
      \imagewithcentreline[width=\textwidth]{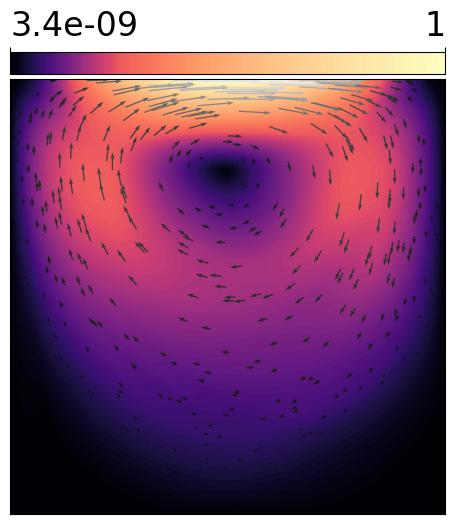}
      \imagewithcentreline[width=\textwidth]{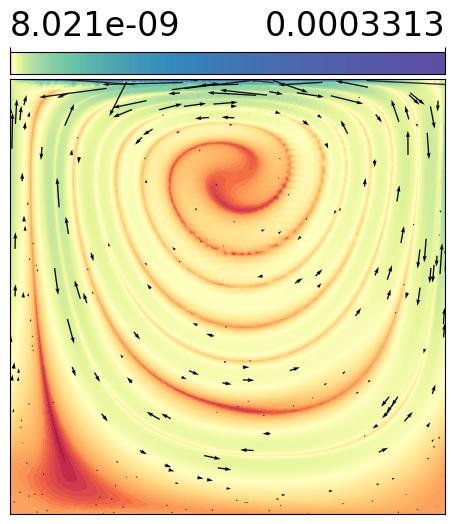} \\
      \imagewithcentreline[width=\textwidth]{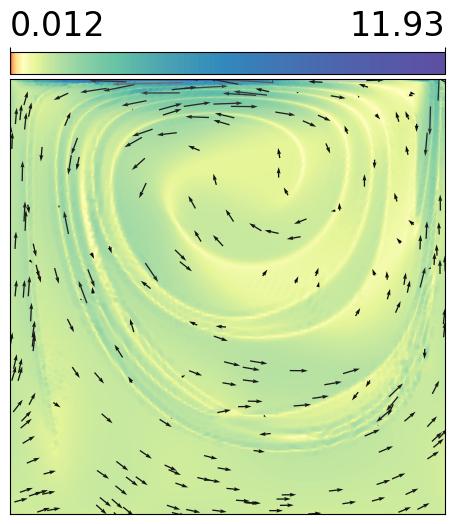}
    \end{subfigure}}
  \subcaptionbox{$t^n = 24.9$}
  {\begin{subfigure}{0.2\textwidth}
      \imagewithcentreline[width=\textwidth]{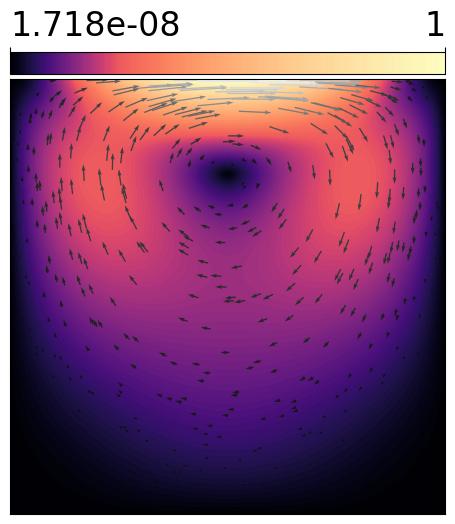} \\
      \imagewithcentreline[width=\textwidth]{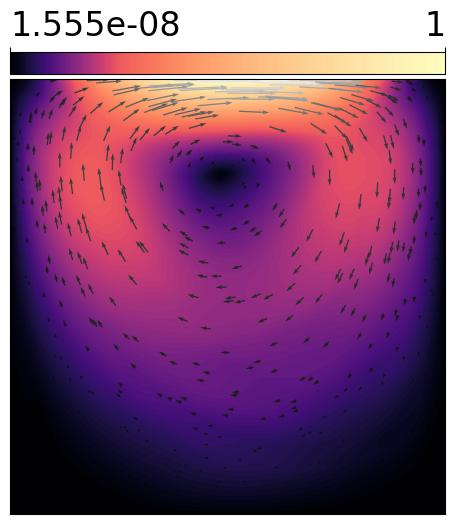}
      \imagewithcentreline[width=\textwidth]{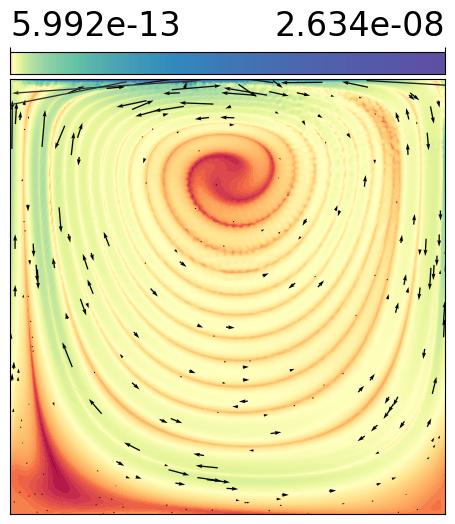} \\
      \imagewithcentreline[width=\textwidth]{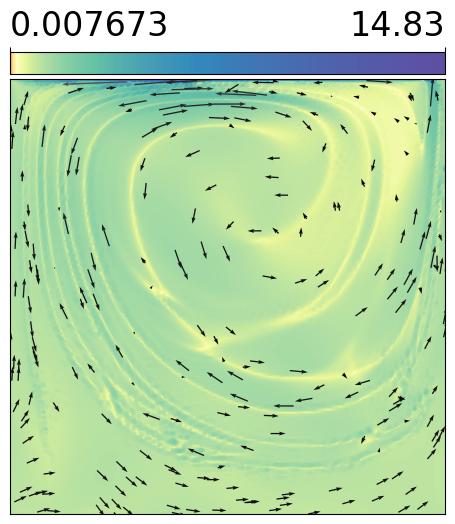}
    \end{subfigure}}
  \caption{Snapshots of $\vec{u}_h$ and $\vec{b}_h$ from the lid-driven cavity experiment in Section \ref{oldroyd-sec:lid-driven-cavity}, with $\Wi = 1$. The top and second rows show $\vec{u}_h$ without and with the extra forcing term, respectively, and the third and bottom rows show $\vec{b}_h$ without and with the extra forcing term, respectively. We draw particular attention to the logarithmic colour scales on each plot, which are used due to the nature of the boundary layer that arises along the top of the domain.}
  \label{oldroyd-fig:lid-driven-cavity-snapshots-large-Wi}
\end{figure}

Along the top moving boundary $\Gamma_D$ we prescribe
$\vec{u} = \Transpose{\qp{u_D, 0}}$, where
\begin{equation}
  u_D\qp{x, t} = 8\qp{1 + \tanh{8\qp{2t-1/2}}}x^2\qp{1 - x}^2,
\end{equation}
and on the remaining boundary $\Gamma_0$ we impose no-slip and
no-penetration conditions $\vec{u} = \vec{0}$. The initial conditions
are $\vec{u}_0 = \vec{0}$ and $\vec{b}_0 = \Transpose{\qp{1, 0}}$.
Representative snapshots of the resulting velocity field and polymeric
vector field (interpreted as the principal stretch direction) for
$\Wi = 1/2$ and $\Wi = 1$ are shown in
Figures~\ref{oldroyd-fig:lid-driven-cavity-snapshots-small-Wi} and
\ref{oldroyd-fig:lid-driven-cavity-snapshots-large-Wi}, respectively.
The centre of rotation for a Newtonian cavity flow is indicated for
reference.

Without the inclusion of the extra forcing term we observe large
spiral-like structures in the profile of $\vec{b}_h$, which can be
interpreted as a visualisation of polymer chains within the fluid.
They decrease in magnitude but become longer and more pronounced as
time evolves. The relaxation term ultimately dominates the stress
generated at the lid and $\vec{b}_h$ is driven towards zero.
Consequently, there is little-to-no visible viscoelastic influence on
$\vec{u}_h$. In particular, we do not observe the upward-left shift of
the centre of fluid rotation (cf.\ Figure~2 in
\cite{ashby2025discretisation}), which is characteristic of
viscoelastic flows in lid-driven cavities. We also highlight the rapid
changes in the direction of $\vec{b}_h$ between neighbouring ``spiral
arms''. At these sharp interfaces it is necessary that
$\abs{\vec{b}_h} < 1$, regardless of the overall magnitude of the
profile, indicating local polymer contraction and thus violating the
assumption of a strongly elongated regime.

With the extra forcing term included, stress generated at the lid is
sustained into the bulk and more complex flow patterns emerge. Similar
rapid transitions in the direction of $\vec{b}_h$ are observed, and
their influence is clearly visible in the velocity field. The centre
of fluid rotation moves upwards and left, although this is a transient
effect and the velocity profile continues to adjust in response to the
unsteady stress. Although the magnitude of $\vec{b}_h$ is generally
slightly larger in the $\Wi = 1$ case (with a maximum of
$\approx 24$ compared to $\approx 21$ when $\Wi = 1/2$), the impact on
the velocity is more pronounced in the smaller Weissenberg regime.
This is consistent with the structure of the coupling term
$(1-\beta)\Wi^{-1}(\vec b\cdot\nabla)\vec b$ in
\eqref{eq:reduced_model_divfree}: both the prefactor $\Wi^{-1}$ and
the spatial gradients of $\vec b_h$ determine the effective forcing on
the momentum equation, and in these simulations their combined effect
is stronger at $\Wi = 1/2$ despite the slightly larger peak values of
$\abs{\vec{b}_h}$ at $\Wi = 1$. The recirculating flow may also play a
role: opposing directions in the principal stress vector can promote
local polymer contraction and thereby limit the growth of stress.

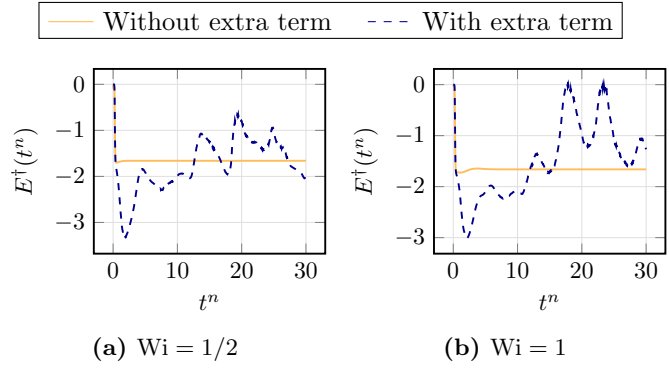
\begin{figure}[h!]
  \centering
  \begin{tikzpicture}
    \node [draw,fill=white] at (0.8,0.3) {\shortstack[l]{
        \ref{no-extra} Without extra term \quad
        \ref{extra} With extra term}};
  \end{tikzpicture}
  \vspace{1em}
  
  \subcaptionbox{$\Wi = 1/2$\label{oldroyd-fig:lid-energy-small}}
  {\begin{tikzpicture}[scale=0.9]
  \begin{axis}[
      cycle list/Dark2,
      thick,
      xlabel=$t^n$,
      ylabel=$E^{\dagger}\qp{t^n}$,
      grid=both,
      minor grid style={gray!25},
      major grid style={gray!25},
      legend style=none,
    ]
    \addplot+[color=small-Wi] table[x=time, y=energy_change, col sep=comma] {data/lid_small_energy.csv};
    \label{no-extra}
    \addplot+[color=large-Wi, dashed] table[x=time, y=energy_change, col sep=comma] {data/lid_small_extra_energy.csv};
    \label{extra}
  \end{axis}
\end{tikzpicture}


}
  \subcaptionbox{$\Wi = 1$\label{oldroyd-fig:lid-energy-large}}
  {\begin{tikzpicture}[scale=0.9]
  \begin{axis}[
      cycle list/Dark2,
      thick,
      xlabel=$t^n$,
      ylabel=$E^{\dagger}\qp{t^n}$,
      grid=both,
      minor grid style={gray!25},
      major grid style={gray!25},
      legend style=none,
    ]
    \addplot+[color=small-Wi] table[x=time, y=energy_change, col sep=comma] {data/lid_large_energy.csv};
    \label{no-extra}
    \addplot+[color=large-Wi, dashed] table[x=time, y=energy_change, col sep=comma] {data/lid_large_extra_energy.csv};
    \label{extra}
  \end{axis}
\end{tikzpicture}


}
  \caption{Plots of the instantaneous rate of change of the modified energy (\ref{oldroyd-eq:modified-energy}) in the lid-driven cavity experiments of Section \ref{oldroyd-sec:lid-driven-cavity}. We emphasise that all values remain negative after initial time, and apparent contradictions are only plotting artifacts.}
  \label{oldroyd-fig:lid-energy-plots}
\end{figure}

In Figure~\ref{oldroyd-fig:lid-energy-plots} we plot the evolution of
the quantity
\begin{equation}\label{oldroyd-eq:modified-energy}
  \begin{aligned}
    E_h^{\dagger}\qp{t^n}
    :=
    E_h\qp{t^n}
    -
    E_h\qp{t^{n-1}}
    +
    \tau\int_{\partial\Omega}
    \big(
    &\Re\vec{u}_h\otimes\vec{u}_h : \vec{u}_h\otimes\normal
      + p_h\normal\cdot\vec{u}_h \\
    &- 2\beta\symgrad{\vec{u}_h} : \vec{u}_h\otimes\normal
      + \vec{u}_h\otimes\vec{b}_h : \vec{b}_h\otimes\normal
      \big) \mathrm d s,
  \end{aligned}
\end{equation}
which represents the instantaneous change of a modified energy that
includes boundary contributions (a quantity that will also be relevant
in the subsequent numerical experiments). The decay of the solution
without the extra term towards a steady state, and the non-steady
behaviour of the solution with the extra term, are both visible in the
profiles of $E_h^{\dagger}$. In all cases $E_h^{\dagger}$ remains
non-positive, demonstrating energy decay consistent with an
appropriately modified version of
Theorem~\ref{thm:discrete-well-posedness-energy}.

\subsection{Flow through a pipe with a cavity}
\label{oldroyd-sec:cavity-flow}

The second physically-motivated experiment is a natural extension of
the lid-driven cavity, featuring flow through a pipe with a large
cavity along one side. Modelling this situation is of practical
interest for investigating the effects of elastic turbulence, which
could be used to clear debris or other unwanted particles from the
cavity region if introduced in an appropriate manner. Reentrant
corners present a challenge to stability through flow separation and
the influence on $\vec{b}$ of the gradient of the fluid velocity,
which can become very large in the vicinity of the reentrant
corner. Analytical results indicate that the velocity vanishes and the
stress becomes singular at these points \cite{hinch1993flow}.

\begin{figure}[h]
  \centering
  \begin{tikzpicture}
  \def\length{10.0};
  \def\height{2.0};
  \def\cavityposition{0.5};
  \def\cavitylength{3.0};
  \def\cavityheight{2.0};

  \colorlet{inflowcolour}{blue}
  \colorlet{outflowcolour}{red}
  \colorlet{noslipcolour}{black}
  
  \begin{scope}[>=Latex]
    \coordinate (inflow-corner-bottom) at (0.0, 0.0);
    \draw[inflowcolour, very thick] (inflow-corner-bottom) -- ([shift={(0.0, \height)}]inflow-corner-bottom) coordinate (inflow-corner-top) node[black, midway, left] {$\inflow$};
    \draw[noslipcolour] (inflow-corner-top) -- ([shift={(\length, 0.0)}]inflow-corner-top) coordinate (outflow-corner-top);
    \draw[outflowcolour, dashed, very thick] (outflow-corner-top) -- ([shift={(0.0, -\height)}]outflow-corner-top) coordinate (outflow-corner-bottom) node[black, midway, right] {$\outflow$};
    \draw[noslipcolour] (outflow-corner-bottom) -- ({\cavityposition * \length + 0.5 * \cavitylength}, 0.0) coordinate (reflex-corner-right);
    \draw[noslipcolour] (reflex-corner-right) -- ([shift={(0.0, -\cavityheight)}]reflex-corner-right) coordinate (cavity-corner-bottom-right);
    \draw[noslipcolour] (cavity-corner-bottom-right) -- ([shift={(-\cavitylength, 0.0)}]cavity-corner-bottom-right) coordinate (cavity-corner-bottom-left);
    \draw[noslipcolour] (cavity-corner-bottom-left) -- ([shift={(0.0, \cavityheight)}]cavity-corner-bottom-left) coordinate (reflex-corner-left);
    \draw[noslipcolour] (reflex-corner-left) -- (inflow-corner-bottom);

    \fill[inflowcolour, opacity=0.2] (inflow-corner-bottom) rectangle ([shift={(-0.1, 0)}]inflow-corner-top);
    \fill[outflowcolour, opacity=0.2] (outflow-corner-bottom) rectangle ([shift={(0.1, 0)}]outflow-corner-top);
    
    \draw[scale=0.5, rotate=90, domain=-2:2, variable=\y]  plot ({\y + 2.0}, {0.5 * \y*\y - 2.0});
    \draw[-latex] ([shift={(0.0, {0.5 * \height})}]inflow-corner-bottom) -- ([shift={(1.0, {0.5 * \height})}]inflow-corner-bottom) node[right] {$\inflowvelocity$};
    \draw[-latex] ([shift={(0.0, {0.4 * \height})}]inflow-corner-bottom) -- ([shift={(0.96, {0.4 * \height})}]inflow-corner-bottom);
    \draw[-latex] ([shift={(0.0, {0.3 * \height})}]inflow-corner-bottom) -- ([shift={(0.84, {0.3 * \height})}]inflow-corner-bottom);
    \draw[-latex] ([shift={(0.0, {0.2 * \height})}]inflow-corner-bottom) -- ([shift={(0.64, {0.2 * \height})}]inflow-corner-bottom);
    \draw[-latex] ([shift={(0.0, {0.1 * \height})}]inflow-corner-bottom) -- ([shift={(0.36, {0.1 * \height})}]inflow-corner-bottom);
    \draw[-latex] ([shift={(0.0, {0.6 * \height})}]inflow-corner-bottom) -- ([shift={(0.96, {0.6 * \height})}]inflow-corner-bottom);
    \draw[-latex] ([shift={(0.0, {0.7 * \height})}]inflow-corner-bottom) -- ([shift={(0.84, {0.7 * \height})}]inflow-corner-bottom);
    \draw[-latex] ([shift={(0.0, {0.8 * \height})}]inflow-corner-bottom) -- ([shift={(0.64, {0.8 * \height})}]inflow-corner-bottom);
    \draw[-latex] ([shift={(0.0, {0.9 * \height})}]inflow-corner-bottom) -- ([shift={(0.36, {0.9 * \height})}]inflow-corner-bottom);
    \draw[-latex] ([shift={({-0.1 * \length}, {0.5 * \height})}]outflow-corner-bottom) node[left] {$\outflowvelocity$} -- ([shift={(0.0, {0.5 * \height})}]outflow-corner-bottom);

    \node at ([shift={(-0.4, -0.4)}]reflex-corner-left) {$\Gamma_0$};
    
    \draw[|<->|] ([shift={(0.0, 0.2)}]inflow-corner-top) -- ([shift={(0.0, 0.2)}]outflow-corner-top) node[midway, above] {$L_1$};
    \draw[|<->|] ([shift={(0.0, -0.2)}]cavity-corner-bottom-left) -- ([shift={(0.0, -0.2)}]cavity-corner-bottom-right) node[midway, below] {$L_2$};
    \draw[|<->|] ([shift={(0.2, 0.)}]cavity-corner-bottom-right) -- ([shift={(0.2, 0.0)}]reflex-corner-right) node[midway, right] {$H_2$};
    \draw[|<->|] ([shift={(0.2, 0.)}]reflex-corner-right) -- ([shift={(0.2, \height)}]reflex-corner-right) node[midway, right] {$H_1$};
    
  \draw plot[smooth cycle, -latex, tension=0.9] coordinates {(0.5 * \length, -0.1 * \cavityheight) (0.5 * \length + 0.35 * \cavitylength, -0.3 * \cavityheight) (0.5 * \length, -0.8 * \cavityheight) (0.5 * \length - 0.35 * \cavitylength, -0.3 * \cavityheight) } [arrow inside={end=latex, opt={scale=1.5}} {0.15, 0.5, 0.8}];
  \draw plot[smooth cycle, -latex, tension=1.0] coordinates {(0.5 * \length, -0.18 * \cavityheight) (0.5 * \length + 0.25 * \cavitylength, -0.35 * \cavityheight) (0.5 * \length, -0.6 * \cavityheight) (0.5 * \length - 0.25 * \cavitylength, -0.35 * \cavityheight) } [arrow inside={end=latex, opt={scale=1.5}} {0.4, 0.9}];
  \draw plot[smooth cycle, -latex, tension=1.0] coordinates {(0.5 * \length, -0.25 * \height) (0.5 * \length + 0.15 * \cavitylength, -0.4 * \cavityheight) (0.5 * \length, -0.5 * \height) (0.5 * \length - 0.15 * \cavitylength, -0.4 * \cavityheight)} [arrow inside={end=latex, opt={scale=1.5}} {0.3}];
  \end{scope}
\end{tikzpicture}
  \caption{A sketch of the cavity flow problem setup in Section \ref{oldroyd-sec:cavity-flow}. The blue line represents the fluid inflow $\inflow$, and the red dashed line represents the outflow $\outflow$. The remaining boundary $\Gamma_0$ represents fixed walls.}
  \label{oldroyd-fig:cavity-setup}
\end{figure}
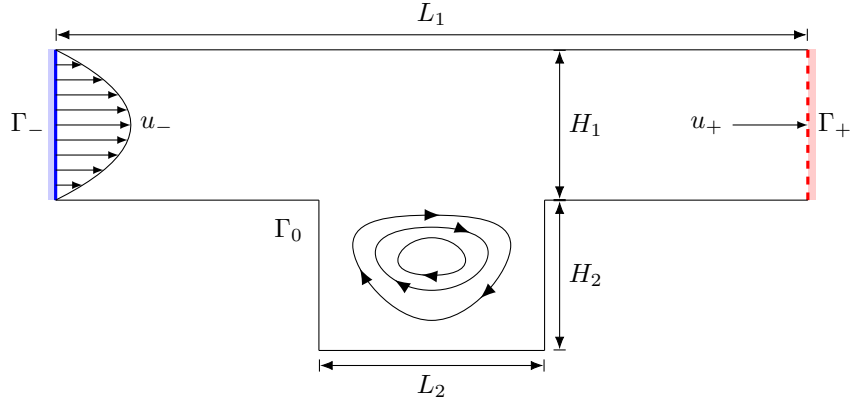

A sketch of the problem setup is shown in Figure
\ref{oldroyd-fig:cavity-setup}, where we take $H_1 = H_2 = 1$, $L_1 =
6$, and $L_2 = 1.5$. The length of the pipe $L_1$ is taken greater
than the other dimensions to allow the flow to fully develop before
reaching the cavity. Taking $T = 20$ and $\tau = 5\times10^{-2}$, we
solve \eqref{eq:fulldis} using polynomial degree $k = 1$ on a
triangulation of the problem domain with increased resolution close to
the boundary (139,341 total degrees of freedom). The mesh is
particularly refined ($h_K \approx 2\times10^{-3}$) around the
reentrant corners either side of the cavity, and the maximum mesh size
is $h_K \approx 6\times10^{-2}$. We take $\sigma_{\mathrm{i}} = 150$
and $\sigma_{\mathrm{b}} = 10^6$. To simulate a creeping flow of a
viscoelastic fluid, the problem parameter values are fixed as $\Re =
5\times10^{-2}$, $\Wi = 1$, and we examine three polymer-to-solvent
ratios: $\beta = 0.9$, $\beta = 0.5$, $\beta = 0.1$. The velocity is
prescribed as parabolic on the inflow $\Gamma_-$ as $\vec{u} =
\Transpose{\qp{u_-, 0}}$, with
\begin{equation}
  \label{oldroyd-eq:inflow-vel-def}
  u_-\qp{y, t} = \frac{4}{H_1^2}\tanh{40t}y\qp{H_1 - y}.
\end{equation}
On the outflow $\Gamma_+$, the natural boundary condition
$-\nabla\vec{u}\cdot\normal + p\normal = \vec{0}$ reduces to the
do-nothing condition $\nabla\vec{u}\cdot\normal = \vec{0}$ by fixing
$p = 0$. On the remaining boundaries, denoted $\Gamma_0$, we set the
no-slip, no-penetration condition $\vec{u} = \vec{0}$. We set $\vec{b}
= \Transpose{\qp{1, 0}}$ on the inflow boundary. The initial
conditions are the same as in Section
\ref{oldroyd-sec:lid-driven-cavity}, with $\vec{u}_0 = \vec{0}$ and
$\vec{b}_0 = \Transpose{\qp{1, 0}}$. We note that, due to the
do-nothing condition on the outflow, the pressure is fully-determined
so no mean-zero condition is required.

\begin{figure}[h!]
  \centering
  \subcaptionbox{$\vec{b}_h$}
  {\begin{subfigure}{0.49\textwidth}
      \includegraphics[width=\textwidth]{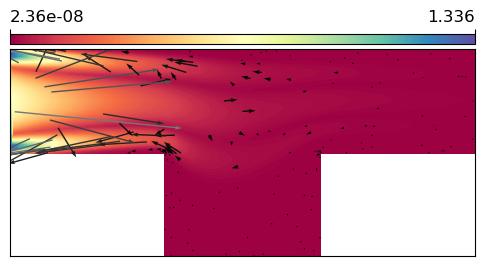} \\
      \includegraphics[width=\textwidth]{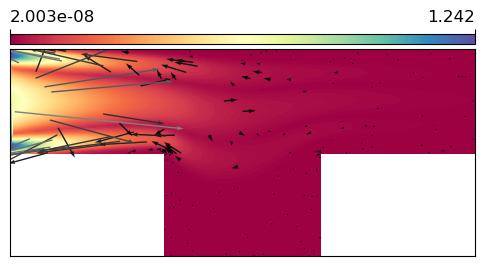} \\
      \includegraphics[width=\textwidth]{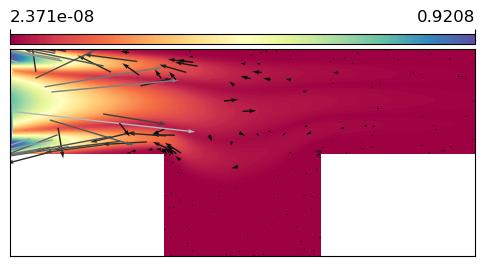}
    \end{subfigure}}
  \subcaptionbox{$\vec{u}_h$}
  {\begin{subfigure}{0.49\textwidth}
      \includegraphics[width=\textwidth]{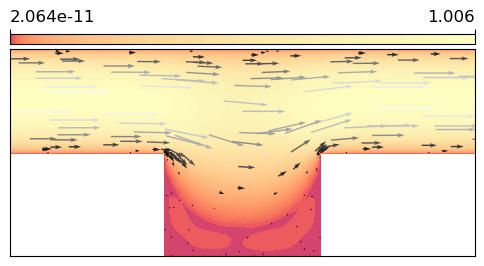} \\
      \includegraphics[width=\textwidth]{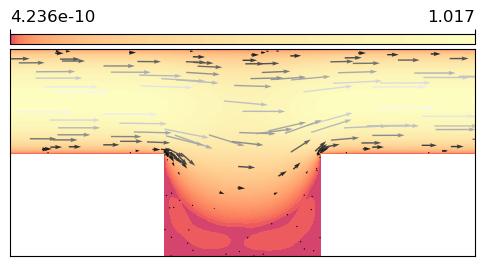} \\
      \includegraphics[width=\textwidth]{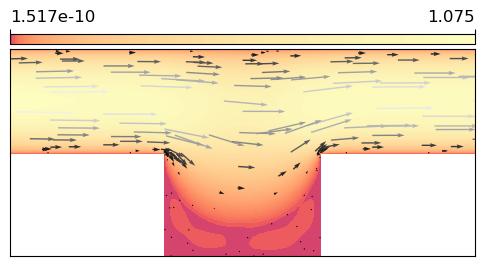}
    \end{subfigure}}
  \caption{Snapshots at the final time of $\vec{b}_h$ and $\vec{u}_h$ with (top) $\beta = 0.9$, (middle) $\beta = 0.5$, and (bottom) $\beta = 0.1$ from the pipe-with-cavity experiments in Section \ref{oldroyd-sec:cavity-flow} without the extra term included. The velocity plots are coloured on a logarithmic scale to emphasise the flow structures present in the cavity.}
  \label{oldroyd-fig:cavity-flow-snapshots}
\end{figure}

\begin{figure}[h!]
  \centering
  \subcaptionbox{$\vec{b}_h$}
  {\begin{subfigure}{0.49\textwidth}
      \includegraphics[width=\textwidth]{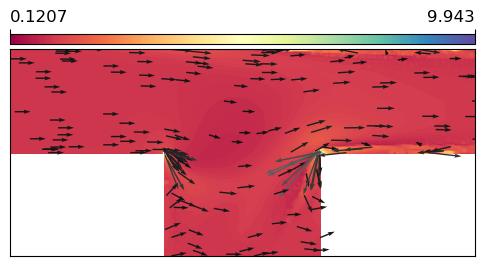} \\
      \includegraphics[width=\textwidth]{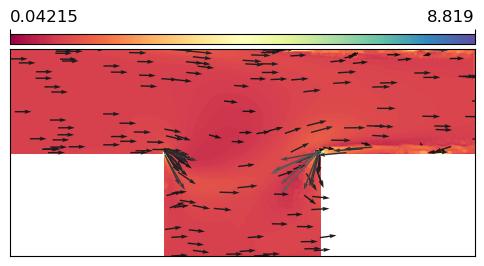} \\
      \includegraphics[width=\textwidth]{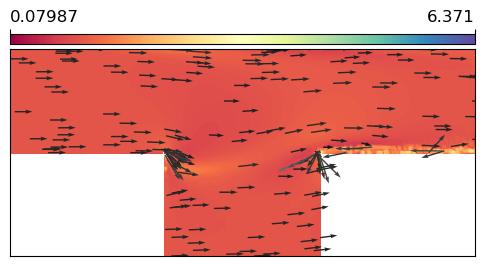}
    \end{subfigure}}
  \subcaptionbox{$\vec{u}_h$}
  {\begin{subfigure}{0.49\textwidth}
      \includegraphics[width=\textwidth]{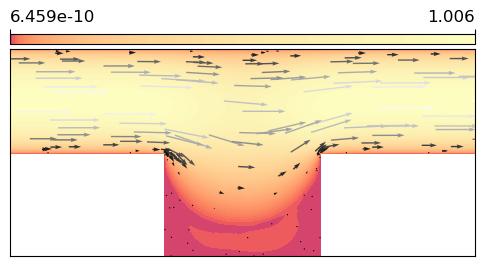} \\
      \includegraphics[width=\textwidth]{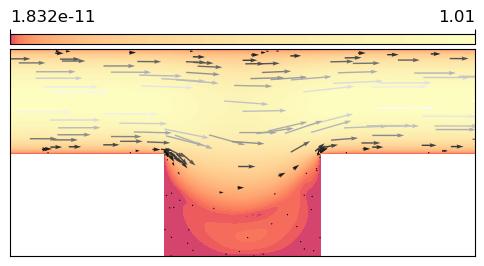} \\
      \includegraphics[width=\textwidth]{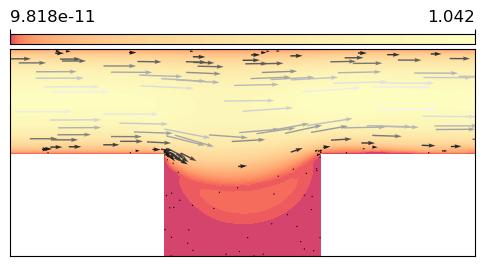}
    \end{subfigure}}
  \caption{Snapshots at the final time of $\vec{b}_h$ and $\vec{u}_h$ with (top) $\beta = 0.9$, (middle) $\beta = 0.5$, and (bottom) $\beta = 0.1$ from the pipe-with-cavity experiments in Section \ref{oldroyd-sec:cavity-flow} with the extra term included. The velocity plots are coloured on a logarithmic scale to emphasise the flow structures present in the cavity.}
  \label{oldroyd-fig:cavity-flow-snapshots-extra}
\end{figure}

The resulting polymeric stress and velocity profiles are depicted at final time without the extra term in Figure \ref{oldroyd-fig:cavity-flow-snapshots} and with the extra term included in Figure \ref{oldroyd-fig:cavity-flow-snapshots-extra}. In both cases it can be seen that the maximum stress attained \emph{decreases} as $\beta$ does, whereas the maximum velocity increases and becomes more affected by the viscoelastic properties of the fluid. Without the extra term present, the stress set at the inflow quickly dies away towards zero as the flow progresses into the domain. Consequently, the velocity profile in these cases remains very similar to that of a Newtonian fluid (e.g. symmetric vortices in the cavity), with only minor deviations visible in the $\beta = 0.1$ case. On the other hand, when the extra term is included the stress is preserved into the bulk and the reentrant corners at the top of the cavity generate large local stress. This feeds back into the velocity profile, and it can be seen that as $\beta$ decreases, the overall velocity in the cavity decreases and a single centre of fluid rotation emerges, breaking the previously-observed symmetry.

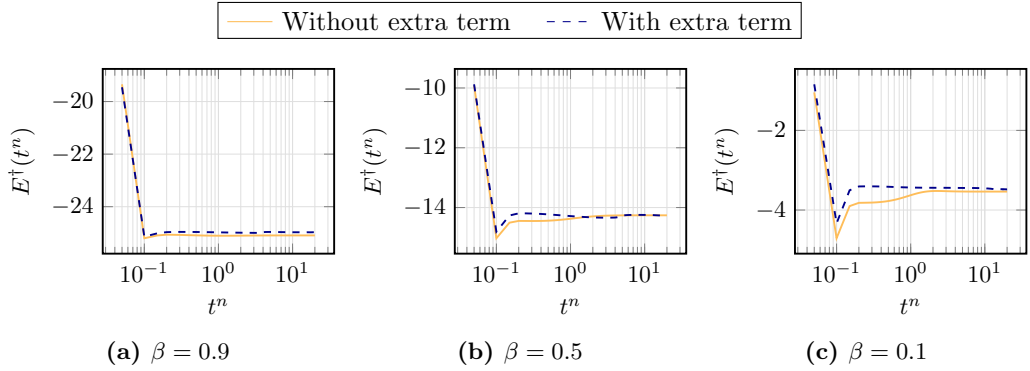
\begin{figure}[h!]
  \centering
  \begin{tikzpicture}
    \node [draw,fill=white] at (0.8,0.3) {\shortstack[l]{
        \ref{no-extra} Without extra term \quad
        \ref{extra} With extra term}};
  \end{tikzpicture}
  \vspace{1em}
  
  \subcaptionbox{$\beta = 0.9$\label{oldroyd-fig:cavity-energy-small}}
  {\begin{tikzpicture}[scale=0.9]
  \begin{axis}[
      cycle list/Dark2,
      thick,
      xmode=log,
      xlabel=$t^n$,
      ylabel=$E^{\dagger}\qp{t^n}$,
      grid=both,
      minor grid style={gray!25},
      major grid style={gray!25},
      legend style=none,
    ]
    \addplot+[color=small-Wi] table[x=time, y=energy_change, col sep=comma] {data/cavity_small_energy.csv};
    \label{no-extra}
    \addplot+[color=large-Wi, dashed] table[x=time, y=energy_change, col sep=comma] {data/cavity_small_extra_energy.csv};
    \label{extra}
  \end{axis}
\end{tikzpicture}


}
  \subcaptionbox{$\beta = 0.5$\label{oldroyd-fig:cavity-energy-moderate}}
  {\begin{tikzpicture}[scale=0.9]
  \begin{axis}[
      cycle list/Dark2,
      thick,
      xmode=log,
      xlabel=$t^n$,
      ylabel=$E^{\dagger}\qp{t^n}$,
      grid=both,
      minor grid style={gray!25},
      major grid style={gray!25},
      legend style=none,
    ]
    \addplot+[color=small-Wi] table[x=time, y=energy_change, col sep=comma] {data/cavity_moderate_energy.csv};
    \label{no-extra}
    \addplot+[color=large-Wi, dashed] table[x=time, y=energy_change, col sep=comma] {data/cavity_moderate_extra_energy.csv};
    \label{extra}
  \end{axis}
\end{tikzpicture}


}
  \subcaptionbox{$\beta = 0.1$\label{oldroyd-fig:cavity-energy-large}}
  {\begin{tikzpicture}[scale=0.9]
  \begin{axis}[
      cycle list/Dark2,
      thick,
      xmode=log,
      xlabel=$t^n$,
      ylabel=$E^{\dagger}\qp{t^n}$,
      grid=both,
      minor grid style={gray!25},
      major grid style={gray!25},
      legend style=none,
    ]
    \addplot+[color=small-Wi] table[x=time, y=energy_change, col sep=comma] {data/cavity_large_energy.csv};
    \label{no-extra}
    \addplot+[color=large-Wi, dashed] table[x=time, y=energy_change, col sep=comma] {data/cavity_large_extra_energy.csv};
    \label{extra}
  \end{axis}
\end{tikzpicture}


}
  \caption{Plots of the instantaneous rate of change of the modified energy (\ref{oldroyd-eq:modified-energy}) in the pipe-with-cavity experiments of Section \ref{oldroyd-sec:cavity-flow}.}
  \label{oldroyd-fig:cavity-energy-plots}
\end{figure}

We plot the instantaneous energy change for this problem in Figure \ref{oldroyd-fig:cavity-energy-plots}. The magnitude of the energy change decreases with $\beta$, since the velocity dissipation term loses significance, emphasising the viscoelasticity of the fluid. In all cases the energy is monotonically decreasing.

\subsection{4:1 planar contraction flows}
\label{oldroyd-sec:contraction-flow}

Another class of problems typically used to benchmark numerical methods for viscoelastic flows are planar contraction flows, in which a fluid flows through a two-dimensional channel featuring a sudden width constriction. Physical experiments highlight the wide range of flow dynamics that can be observed by slight adjustments to the geometry or fluid properties \cite{evans1986flow}.

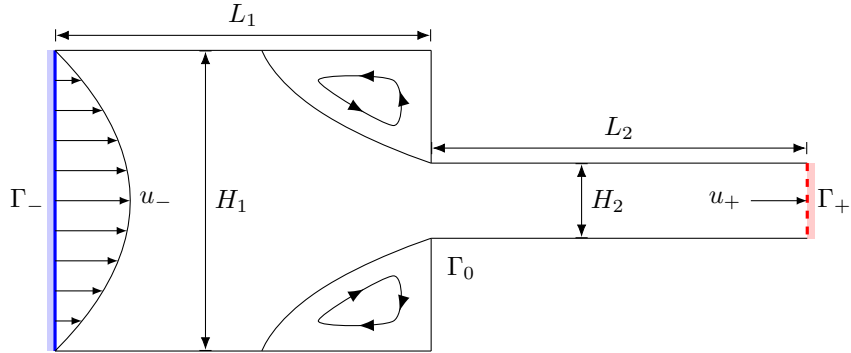
\begin{figure}[h]
  \centering
  \begin{tikzpicture}
  \def\varupstreamlength{5.0};
  \def\varheight{4.0};
  \def\varcontractionposition{0.5};
  \def\varcontractionlength{5.0};
  \def\varcontractionheight{1.0};

  \colorlet{inflowcolour}{blue}
  \colorlet{outflowcolour}{red}
  \colorlet{noslipcolour}{black}
  
  \begin{scope}[>=Latex]
    \coordinate (inflow-corner-bottom) at (0.0, 0.0);
    \draw[inflowcolour, very thick] (inflow-corner-bottom) -- ([shift={(0.0, \varheight)}]inflow-corner-bottom) coordinate (inflow-corner-top) node[black, midway, left] {$\inflow$};
    \draw[noslipcolour] (inflow-corner-top) -- ([shift={(\varupstreamlength, 0.0)}]inflow-corner-top) coordinate (upstream-corner-top);
    \draw[noslipcolour] (upstream-corner-top) -- ([shift={(0.0, {-0.5 * \varheight + \varcontractionposition * \varcontractionheight})}]upstream-corner-top) coordinate (upstream-contraction-corner-top);
    \draw[noslipcolour] (upstream-contraction-corner-top) -- ([shift={(\varcontractionlength, 0.0)}]upstream-contraction-corner-top) coordinate (outflow-corner-top);
    \draw[outflowcolour, dashed, very thick] (outflow-corner-top) -- ([shift={(0.0, -\varcontractionheight)}]outflow-corner-top) coordinate (outflow-corner-bottom) node[black, midway, right] {$\outflow$};
    \draw[noslipcolour] (outflow-corner-bottom) -- ([shift={(-\varcontractionlength, 0.0)}]outflow-corner-bottom) coordinate (upstream-contraction-corner-bottom);
    \draw[noslipcolour] (upstream-contraction-corner-bottom) -- ([shift={(0.0, {-0.5 * \varheight + \varcontractionposition * \varcontractionheight})}]upstream-contraction-corner-bottom) coordinate (upstream-corner-bottom);
    \draw[noslipcolour] (upstream-corner-bottom) -- ([shift={(-\varupstreamlength, 0.0)}]upstream-corner-bottom) coordinate (inflow-corner-bottom);
    
    \fill[inflowcolour, opacity=0.2] (inflow-corner-bottom) rectangle ([shift={(-0.1, 0)}]inflow-corner-top);
    \fill[outflowcolour, opacity=0.2] (outflow-corner-bottom) rectangle ([shift={(0.1, 0)}]outflow-corner-top);
    
    \draw[|<->|] ([shift={(0.0, 0.2)}]inflow-corner-top) -- ([shift={(0.0, 0.2)}]upstream-corner-top) node[midway, above] {$\upstreamlength$};
    \draw[|<->|] ([shift={(0.0, 0.2)}]upstream-contraction-corner-top) -- ([shift={(0.0, 0.2)}]outflow-corner-top) node[midway, above] {$\contractionlength$};
    \draw[|<->|] ([shift={({0.4 * \varupstreamlength}, 0.0)}]inflow-corner-top) -- ([shift={({0.4 * \varupstreamlength}, 0.0)}]inflow-corner-bottom) node[midway, right] {$\upstreamwidth$};
    \draw[|<->|] ([shift={({-0.6 * \varcontractionlength}, 0.0)}]outflow-corner-top) -- ([shift={({-0.6 * \varcontractionlength}, 0.0)}]outflow-corner-bottom) node[midway, right] {$\contractionwidth$};

    \draw[scale=1.0, rotate=90, domain=-2:2, variable=\y]  plot ({\y + 2.0}, {0.25 * \y*\y - 1.0});
    \draw[-latex] ([shift={(0.0, {0.5 * \varheight})}]inflow-corner-bottom) -- ([shift={(1.0, {0.5 * \varheight})}]inflow-corner-bottom) node[right] {$\inflowvelocity$};
    \draw[-latex] ([shift={(0.0, {0.4 * \varheight})}]inflow-corner-bottom) -- ([shift={(0.96, {0.4 * \varheight})}]inflow-corner-bottom);
    \draw[-latex] ([shift={(0.0, {0.3 * \varheight})}]inflow-corner-bottom) -- ([shift={(0.84, {0.3 * \varheight})}]inflow-corner-bottom);
    \draw[-latex] ([shift={(0.0, {0.2 * \varheight})}]inflow-corner-bottom) -- ([shift={(0.64, {0.2 * \varheight})}]inflow-corner-bottom);
    \draw[-latex] ([shift={(0.0, {0.1 * \varheight})}]inflow-corner-bottom) -- ([shift={(0.36, {0.1 * \varheight})}]inflow-corner-bottom);
    \draw[-latex] ([shift={(0.0, {0.6 * \varheight})}]inflow-corner-bottom) -- ([shift={(0.96, {0.6 * \varheight})}]inflow-corner-bottom);
    \draw[-latex] ([shift={(0.0, {0.7 * \varheight})}]inflow-corner-bottom) -- ([shift={(0.84, {0.7 * \varheight})}]inflow-corner-bottom);
    \draw[-latex] ([shift={(0.0, {0.8 * \varheight})}]inflow-corner-bottom) -- ([shift={(0.64, {0.8 * \varheight})}]inflow-corner-bottom);
    \draw[-latex] ([shift={(0.0, {0.9 * \varheight})}]inflow-corner-bottom) -- ([shift={(0.36, {0.9 * \varheight})}]inflow-corner-bottom);
    \draw[-latex] ([shift={({-0.15 * \varcontractionlength}, {0.5 * \varcontractionheight})}]outflow-corner-bottom) node[left] {$\outflowvelocity$} -- ([shift={(0.0, {0.5 * \varcontractionheight})}]outflow-corner-bottom);

    \node at ([shift={(0.4, -0.4)}]upstream-contraction-corner-bottom) {$\Gamma_0$};
    
    \draw plot[smooth, tension=1] coordinates {([shift={(-0.45 * \varupstreamlength, 0.0)}]upstream-corner-top) ([shift={(-0.3 * \varupstreamlength, -0.2 * \varheight)}]upstream-corner-top) (upstream-contraction-corner-top)};
    \draw plot[smooth cycle, -latex] coordinates {([shift={(-0.3 * \varupstreamlength, -0.1 * \varheight)}]upstream-corner-top) ([shift={(-0.1 * \varupstreamlength, -0.25 * \varheight)}]upstream-corner-top) ([shift={(-0.1 * \varupstreamlength, -0.1 * \varheight)}]upstream-corner-top)} [arrow inside={end=latex, opt={scale=1.5}} {0.18, 0.4, 0.85}];
    \draw plot[smooth, tension=1] coordinates {([shift={(-0.45 * \varupstreamlength, 0.0)}]upstream-corner-bottom) ([shift={(-0.3 * \varupstreamlength, 0.2 * \varheight)}]upstream-corner-bottom) (upstream-contraction-corner-bottom)};
    \draw plot[smooth cycle, -latex] coordinates {([shift={(-0.3 * \varupstreamlength, 0.1 * \varheight)}]upstream-corner-bottom) ([shift={(-0.1 * \varupstreamlength, 0.25 * \varheight)}]upstream-corner-bottom) ([shift={(-0.1 * \varupstreamlength, 0.1 * \varheight)}]upstream-corner-bottom)} [arrow inside={end=latex, opt={scale=1.5}} {0.18, 0.4, 0.85}];
  \end{scope}
\end{tikzpicture}

  \caption{A sketch of the 4:1 contraction flow problem setup in Section \ref{oldroyd-sec:contraction-flow}. The blue line represents the fluid inflow $\inflow$, and the red dashed line represents the outflow $\outflow$. The remaining boundary $\Gamma_0$ represents fixed walls.}
  \label{oldroyd-fig:contraction-setup}
\end{figure}

A sketch of the problem setup can be seen in Figure \ref{oldroyd-fig:contraction-setup}, where we take a pre-contraction channel length of $L_1 = 2$ and width of $\upstreamwidth = 1$ and a post-contraction channel length of $L_2 = 1$ and width of $\contractionwidth = 0.25$. This gives a 4:1 contration ratio ($\upstreamwidth / \contractionwidth$), which is a common choice in the literature. However, the ratio can have important influence, and some flow features are only observed with certain setups \cite{alves2004effect}.

Fixing the final time as $T = 20$ and the time step size as $\tau = 10^{-3}$, we solve \eqref{eq:fulldis}, with polynomial degree $k = 1$, on a triangulation (56,876 total degrees of freedom) of the spatial domain. The triangulation features increased resolution around the boundaries and especially the reentrant corners, which we regularise by rounding to have radius 0.02 and about which $h_K \approx 10^{-3}$. We take $\sigma_{\mathrm{i}} = 300$ and $\sigma_{\mathrm{b}} =10^6$. Taking $\Re = 10^{-2}$ and $\beta = 0.9$, we consider three Weissenberg number regimes: $\Wi = 1.4$, $\Wi = 1.5$, $\Wi = 1.6$. As in the previous pipe-with-cavity experiment, the velocity is prescribed on the inflow boundary $\inflow$ as $\vec{u} = \Transpose{\qp{\inflowvelocity, 0}}$, where $u_-$ is defined in \eqref{oldroyd-eq:inflow-vel-def}, and the do-nothing condition $\nabla\vec{u}\cdot\normal = \vec{0}$ is satisfied on the outflow boundary $\outflow$ by fixing $p = 0$. The remaining boundaries, denoted $\Gamma_0$, feature the no-slip, no-penetration condition $\vec{u} = \vec{0}$. We set $\vec{b} = \Transpose{\qp{2, 0}}$ on the inflow boundary, and the do-nothing condition for $\vec{b}$ is satisfied by fixing $r = 0$ on the remaining boundary. The initial conditions are given as $\vec{u}_0 = \vec{0}$ and $\vec{b}_0 = \Transpose{\qp{2, 0}}$. Again, due to due to setting the pressure on the outflow, it is fully-determined and no mean-zero condition is required.

\begin{figure}[h!]
  \centering
  \subcaptionbox{$\vec{b}_h$}
  {\begin{subfigure}{0.49\textwidth}
      \includegraphics[width=\textwidth]{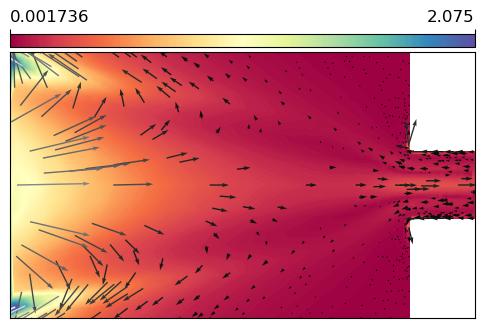} \\
      \includegraphics[width=\textwidth]{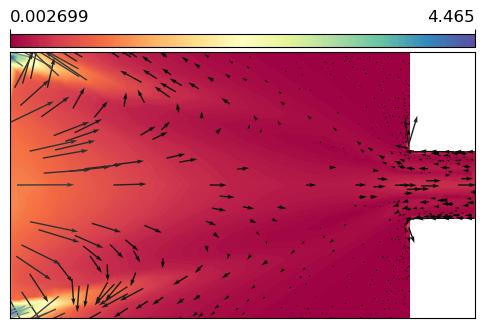} \\
      \includegraphics[width=\textwidth]{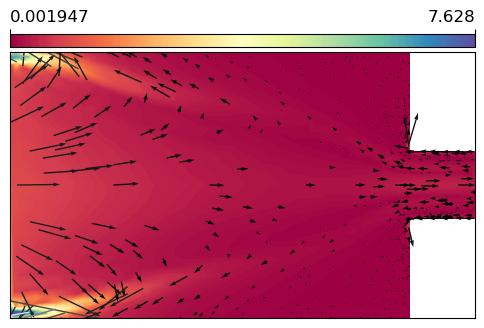}
    \end{subfigure}}
  \subcaptionbox{$\vec{u}_h$}
  {\begin{subfigure}{0.49\textwidth}
      \includegraphics[width=\textwidth]{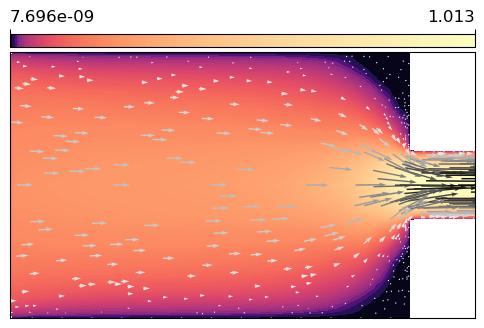} \\
      \includegraphics[width=\textwidth]{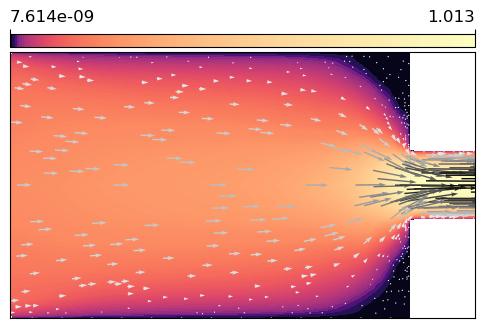} \\
      \includegraphics[width=\textwidth]{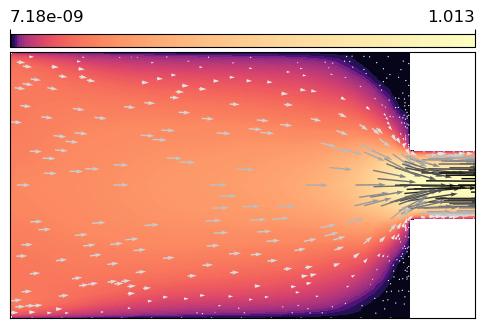}
    \end{subfigure}}
  \caption{Snapshots at $t^n = 5$ of $\vec{b}_h$ and $\vec{u}_h$ with (top) $\Wi = 1.4$, (middle) $\Wi = 1.5$, and (bottom) $\Wi = 1.6$ from the 4:1 contraction experiments in Section \ref{oldroyd-sec:contraction-flow} without the extra term included. The plots are coloured on a logarithmic scale to emphasise the structures present.}
  \label{oldroyd-fig:contraction-flow-snapshots}
\end{figure}

\begin{figure}[h!]
  \centering
  \subcaptionbox{$\vec{b}_h$}
  {\begin{subfigure}{0.49\textwidth}
      \includegraphics[width=\textwidth]{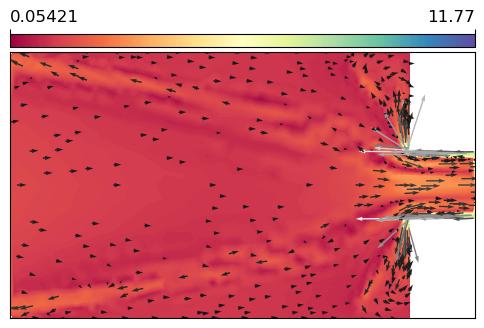} \\
      \includegraphics[width=\textwidth]{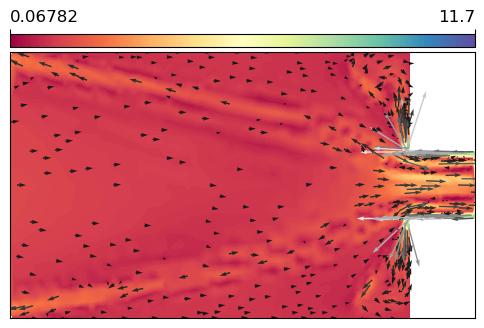} \\
      \includegraphics[width=\textwidth]{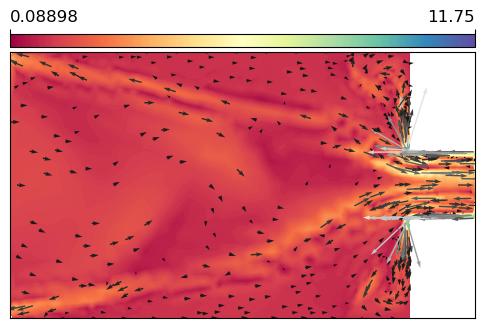}
    \end{subfigure}}
  \subcaptionbox{$\vec{u}_h$}
  {\begin{subfigure}{0.49\textwidth}
      \includegraphics[width=\textwidth]{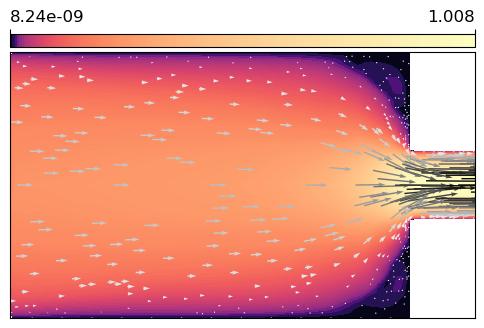} \\
      \includegraphics[width=\textwidth]{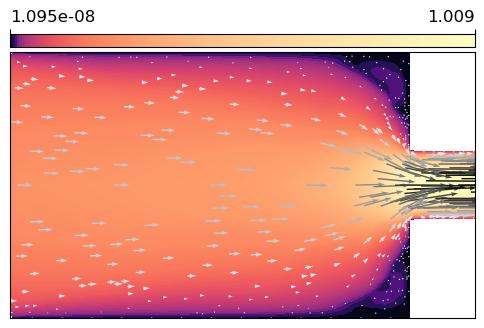} \\
      \includegraphics[width=\textwidth]{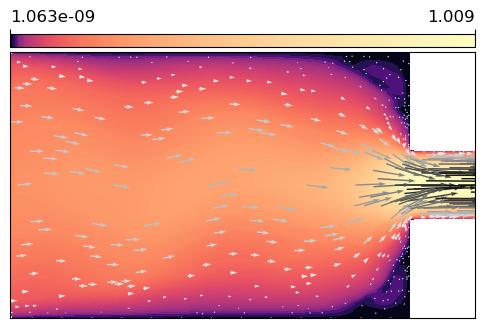}
    \end{subfigure}}
  \caption{Snapshots at $t^n = 5$ of $\vec{b}_h$ and $\vec{u}_h$ with (top) $\Wi = 1.4$, (middle) $\Wi = 1.5$, and (bottom) $\Wi = 1.6$ from the 4:1 contraction experiments in Section \ref{oldroyd-sec:contraction-flow} with the extra term included. The plots are coloured on a logarithmic scale to emphasise the structures present.}
  \label{oldroyd-fig:contraction-flow-snapshots-extra}
\end{figure}

Snapshots of the principal polymeric stress and velocity profiles at $t^n = 5$ are shown without the extra term in Figure \ref{oldroyd-fig:contraction-flow-snapshots} and with the extra term in Figure \ref{oldroyd-fig:contraction-flow-snapshots-extra}. This plotting time is chosen to exemplify the differences between the flows for different Weissenberg numbers. In all cases stress is advected along with the flow and builds up around the reentrant corners. As in the previous numerical experiments, without the extra term the bulk stress decays to zero, but with the extra term two different effects occur after sufficient stress has accumulated. As $\Wi$ increases vortices appear in the top and bottom corners in line with the contraction, where the magnitude of the velocity is small, and they grow in strength. Then, in the $\Wi = 1.6$ case, the flow develops oscillations further upstream of the contraction, which can be observed in both the stress and velocity profiles. Both of these behaviours are in line with results from the literature, with first corner vortices appearing for Weissenberg numbers above a certain value and then oscillations at even higher values.

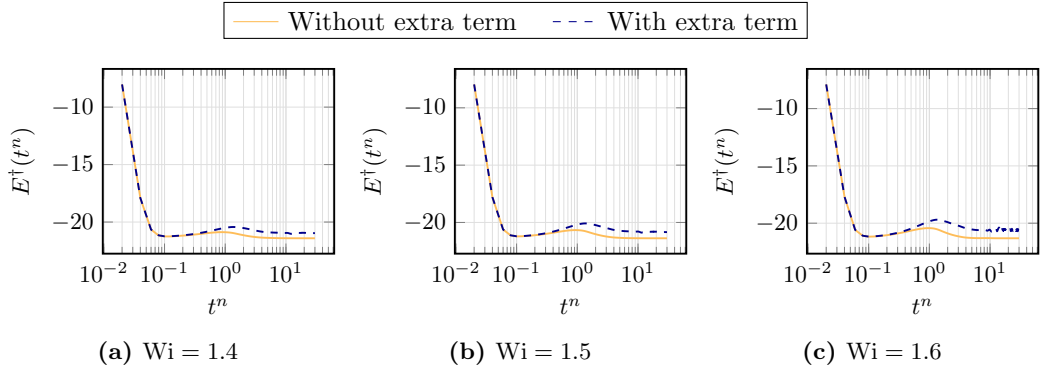
\begin{figure}[h!]
  \centering
  \begin{tikzpicture}
    \node [draw,fill=white] at (0.8,0.3) {\shortstack[l]{
        \ref{no-extra} Without extra term \quad
        \ref{extra} With extra term}};
  \end{tikzpicture}
  \vspace{1em}
  
  \subcaptionbox{$\Wi = 1.4$\label{oldroyd-fig:contraction-energy-small}}
  {\begin{tikzpicture}[scale=0.9]
  \begin{axis}[
      cycle list/Dark2,
      thick,
      xmode=log,
      xlabel=$t^n$,
      ylabel=$E^{\dagger}\qp{t^n}$,
      grid=both,
      minor grid style={gray!25},
      major grid style={gray!25},
      legend style=none,
    ]
    \addplot+[color=small-Wi] table[x=time, y=energy_change, col sep=comma] {data/contraction_small_energy.csv};
    \label{no-extra}
    \addplot+[color=large-Wi, dashed] table[x=time, y=energy_change, col sep=comma] {data/contraction_small_extra_energy.csv};
    \label{extra}
  \end{axis}
\end{tikzpicture}


}
  \subcaptionbox{$\Wi = 1.5$\label{oldroyd-fig:contraction-energy-moderate}}
  {\begin{tikzpicture}[scale=0.9]
  \begin{axis}[
      cycle list/Dark2,
      thick,
      xmode=log,
      xlabel=$t^n$,
      ylabel=$E^{\dagger}\qp{t^n}$,
      grid=both,
      minor grid style={gray!25},
      major grid style={gray!25},
      legend style=none,
    ]
    \addplot+[color=small-Wi] table[x=time, y=energy_change, col sep=comma] {data/contraction_moderate_energy.csv};
    \label{no-extra}
    \addplot+[color=large-Wi, dashed] table[x=time, y=energy_change, col sep=comma] {data/contraction_moderate_extra_energy.csv};
    \label{extra}
  \end{axis}
\end{tikzpicture}


}
  \subcaptionbox{$\Wi = 1.6$\label{oldroyd-fig:contraction-energy-large}}
  {\begin{tikzpicture}[scale=0.9]
  \begin{axis}[
      cycle list/Dark2,
      thick,
      xmode=log,
      xlabel=$t^n$,
      ylabel=$E^{\dagger}\qp{t^n}$,
      grid=both,
      minor grid style={gray!25},
      major grid style={gray!25},
      legend style=none,
    ]
    \addplot+[color=small-Wi] table[x=time, y=energy_change, col sep=comma] {data/contraction_large_energy.csv};
    \label{no-extra}
    \addplot+[color=large-Wi, dashed] table[x=time, y=energy_change, col sep=comma] {data/contraction_large_extra_energy.csv};
    \label{extra}
  \end{axis}
\end{tikzpicture}


}
  \caption{Plots of the instantaneous rate of change of the modified energy (\ref{oldroyd-eq:modified-energy}) in the 4:1 planar contraction flow experiments of Section \ref{oldroyd-sec:contraction-flow}.}
  \label{oldroyd-fig:contraction-energy-plots}
\end{figure}

The instantaneous energy change for the 4:1 contraction flow experiments are shown in Figure \ref{oldroyd-fig:contraction-energy-plots}, where the monotonic energy decrease is demonstrated in all cases. Although the changes are minor between each of the Weissenberg number cases, the impact can be seen as a slight increase as $\Wi$ grows in the experiments including the extra term, indicating that the problem becomes less dissipative.

\section{Conclusions}
\label{sec:conclusion}

In this work, we developed a finite element framework for a uniaxial
reduction of the Oldroyd-B model. We derived spatially semidiscrete
and fully discrete formulations, ensuring compatibility with the
underlying PDE stability framework. Existence and uniqueness of
solutions were shown for the spatially semidiscrete and fully discrete
problems.

Numerical experiments validated the approach on benchmark
problems. These results confirmed the robustness of the method for
handling high Weissenberg numbers and resolving sharp stress
gradients. Despite the simplifications introduced by the uniaxial
reduction, key features of viscoelastic flow dynamics, such as elastic
effects and flow instabilities, were captured. 

Future research could focus on extending the framework to handle
multi-component polymeric flows or incorporating additional physical
effects, such as temperature-dependence or shear-thinning
behaviour. Investigating the limitations of the uniaxial model,
particularly in regimes where isotropic stress components are
significant, is another direction.

\printbibliography

\end{document}